\documentclass [11pt] {article}
\usepackage[T1]{fontenc}
\usepackage{lmodern} 
\usepackage{textcomp}
\usepackage{amsmath}
\usepackage{amsfonts}
\usepackage{amsthm}
\usepackage{amssymb}
\usepackage{graphicx}
\usepackage{multicol}
\usepackage{color}
\usepackage{wrapfig}
\usepackage{subfig}
\usepackage[mathscr]{euscript}
\usepackage[all]{xy}
\usepackage{comment}
\usepackage{enumerate}
\usepackage{manfnt} 
\usepackage{array}
\usepackage{tikz-cd}
\usepackage[hidelinks]{hyperref}
\usepackage{cleveref}
\usepackage{aliascnt}
\usepackage{float}
\usepackage{caption}
\usepackage{microtype}
\usepackage{tabularx}
\usepackage{eso-pic}
\usetikzlibrary{fit, backgrounds}
\usetikzlibrary{calc}
\usetikzlibrary{decorations.markings} 
\usetikzlibrary{arrows.meta}           

\DeclareUnicodeCharacter{221E}{\ensuremath{\infty}}

\newtheoremstyle{nonitalic}
  {3pt}        
  {3pt}        
  {}           
  {}           
  {\bfseries}  
  {.}          
  {.5em}       
  {}           

\theoremstyle{plain}

\newtheorem{theorem}{Theorem}[section]
\crefname{theorem}{Theorem}{Theorems}
\Crefname{theorem}{Theorem}{Theorems}

\crefname{mainthm}{Theorem}{Theorems}
\Crefname{mainthm}{Theorem}{Theorems}

\newaliascnt{proposition}{theorem}
\newtheorem{proposition}[proposition]{Proposition}
\aliascntresetthe{proposition}
\crefname{proposition}{Proposition}{Propositions}
\Crefname{proposition}{Proposition}{Propositions}

\newaliascnt{corollary}{theorem}
\newtheorem{corollary}[corollary]{Corollary}
\aliascntresetthe{corollary}
\crefname{corollary}{Corollary}{Corollaries}
\Crefname{corollary}{Corollary}{Corollaries}

\newaliascnt{lemma}{theorem}
\newtheorem{lemma}[lemma]{Lemma}
\aliascntresetthe{lemma}
\crefname{lemma}{Lemma}{Lemmas}
\Crefname{lemma}{Lemma}{Lemmas}

\newaliascnt{conjecture}{theorem}

\aliascntresetthe{conjecture}
\crefname{conjecture}{Conjecture}{Conjectures}
\Crefname{conjecture}{Conjecture}{Conjectures}

\theoremstyle{nonitalic}

\newaliascnt{definition}{theorem}
\newtheorem{definition}[definition]{Definition}
\aliascntresetthe{definition}
\crefname{definition}{Definition}{Definitions}
\Crefname{definition}{Definition}{Definitions}

\newaliascnt{example}{theorem}
\newtheorem{example}[example]{Example}
\aliascntresetthe{example}
\crefname{example}{Example}{Examples}
\Crefname{example}{Example}{Examples}

\newaliascnt{rem}{theorem}
\newtheorem{rem}[rem]{Remark}
\aliascntresetthe{rem}
\crefname{rem}{Remark}{Remarks}
\Crefname{rem}{Remark}{Remarks}

\newtheorem*{rem*}{Remark}

\newcommand{\R}{\mathbb{R}}

\newcommand{\Q}{\mathbb{Q}}
\newcommand{\C}{\mathbb{C}}
\newcommand{\A}{\mathbb{A}}

\newcommand{\cC}{\mathcal{C}}
\newcommand{\cD}{\mathcal{D}}

\newcommand{\cS}{\mathcal{S}}
\newcommand{\cE}{\mathcal{E}}
\newcommand{\cX}{\mathcal{X}}

\newcommand{\cU}{\mathcal{U}}

\newcommand{\lf}{\left}
\newcommand{\ri}{\right}
\newcommand{\norm}[1]{\left\vert{#1}\right\vert} 
\newcommand{\Norm}[1]{\left\Vert{#1}\right\Vert} 
\newcommand{\Map}{\operatorname{Map}}
\newcommand{\Emb}{\operatorname{Emb}}

\newcommand{\Psh}{\operatorname{PSh}}

\newcommand{\Disk}{\operatorname{Disk}}
\newcommand{\Man}{\operatorname{Mfld}}
\newcommand{\Top}{\operatorname{Top}}
\newcommand{\Strat}{\operatorname{CStrat}}
\newcommand{\cDisk}{\mathcal{D}\mathrm{isk}}
\newcommand{\cMan}{\mathcal{M}\mathrm{fld}}

\newcommand{\Fun}{\operatorname{Fun}}

\newcommand{\Set}{\operatorname{Set}}
\newcommand{\sSet}{\operatorname{sSet}}
\newcommand{\Sing}{\operatorname{Sing}}

\newcommand{\Sd}{\operatorname{Sd}}

\newcommand{\Conf}{\operatorname{Conf}}

\newcommand{\Cat}{\operatorname{Cat}_\infty}

\newcommand{\colim}[2]{\mathop{\operatorname{colim}}\limits_{#1} #2}
\newcommand{\limi}[2]{\mathop{\operatorname{lim}}\limits_{#1} #2}
\newcommand{\Fin}{\operatorname{Fin}}

\newcommand{\holim}[2]{\mathop{\operatorname{holim}}\limits_{#1} #2}

\newcommand{\exit}[2]{\operatorname{Exit}^{#1}(#2)}
\newcommand{\enter}[2]{\operatorname{Enter}^{#1}(#2)}
\newcommand{\cart}{\operatorname{Cart}}
\newcommand{\us}[1]{\underline{#1}}
\newcommand{\la}{\langle}
\newcommand{\ra}{\rangle}
\newcommand{\pos}[1]{\operatorname{B}({#1})}
\newcommand{\prt}[1]{\operatorname{P}({#1})}
\newcommand{\Pos}{\operatorname{Poset}}

\newcommand{\con}{\textbf{con}}
\newcommand{\F}{\mathcal{X}^\leq_I} 
\newcommand{\X}{\mathcal{X}^P} 
\newcommand{\B}{\mathcal{X}^B} 

\tikzset{
  midarrow/.style={
    postaction={decorate},
    decoration={
      markings,
      mark=at position 0.5 with {\arrow{>}}
    }
  }
}

\makeatother

\begin{document}

\title{A simplicial model for manifold calculus through the configuration category}
\date{2026}

\author{Rodrigue Haya Enriquez}


\maketitle

\begin{abstract}
Let $M$ be a closed $d$-dimensional smooth manifold and $X:\Delta^{op}\rightarrow\Fin$ a finite simplicial set whose geometric realization is homeomorphic to $M$. Given the $\infty$-category $\cDisk_{/M}$ of disks embedded in $M$ and the path components functor $\pi_0:\cDisk_{/M}\rightarrow\Fin$, we construct a cofinal functor $\Gamma_X:\Delta^{op}\rightarrow\cDisk_{/M}$ such that $X=\pi_0\circ\Gamma_X$.

This result provides an alternative construction of a simplicial model developed by Boavida de Brito, Lambrechts, Pryor and Songhafouo Tsopméné to compute the space-valued presheaf $\Emb(-,M):\cDisk^{op}\rightarrow\cS$. We also extend the result to manifolds with a $G$-structure and discuss applications to manifold calculus, factorization homology, and a possible generalization to manifolds with boundary.
\end{abstract}

\begingroup
\small
\setlength{\parskip}{2pt}
\linespread{0.9}\selectfont
\tableofcontents
\endgroup


\section{Introduction}

For $d\geq 0$, let $\cMan_d$ be the $\infty$-category of $d$-dimensional smooth manifolds without boundary and embeddings between them. In its modern form initiated in \cite{manifold}, \emph{manifold calculus} studies space-valued presheaves $F:\cMan_d^{op}\rightarrow\cS$ through their restrictions to open balls. Namely, for each $k\geq 1$, consider the full subcategory $\cDisk_{< k}\subset\cMan_d$ whose objects are disjoint unions of less than $k$ $d$-dimensional euclidean open balls. Likewise, let $\cDisk_{<\infty}$ be the full subcategory whose objects are finite disjoint unions of $d$-dimensional euclidean open balls. 

For each presheaf $F:\cMan_d^{op}\rightarrow\cS$, we get a restricted presheaf $F_{|\cDisk_{< k}}:\cDisk_{< k}^{op}\rightarrow\cS$ for all $k\in\mathbb{Z}_{\geq 1}\cup\{+\infty\}$. These restrictions allow to define $T_kF:\cMan_d^{op}\rightarrow\cS$ as a space of natural transformations between $\cDisk_{<k+1}$-presheaves:

\[
T_kF(M):=\Map_{\Psh(\cDisk_{<k+1})}(\Emb(-,M),F).
\]

From the sequence of subcategories
\[
\cDisk_{< 1}\subset\cDisk_{< 2}\subset\cdots\cDisk_{<\infty}\subset\cMan_d,
\]

we then get the \emph{Taylor tower} of $F$

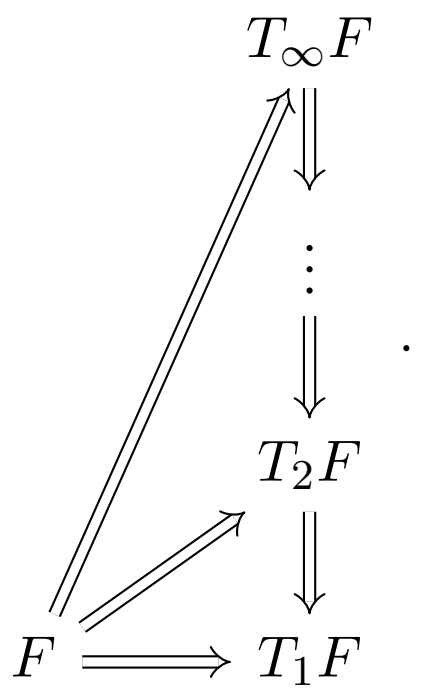
\[
\begin{tikzcd}
& T_\infty F\arrow[d,Rightarrow] \\
& \vdots\arrow[d,Rightarrow] \\
& T_2F\arrow[d,Rightarrow] \\
F\arrow[uuur,Rightarrow]\arrow[r,Rightarrow]\arrow[ur,Rightarrow] & T_1F
\end{tikzcd}.
\]

The benefit of this tower is that in cases of interest, the natural maps $F\Rightarrow T_kF$ become more and more connected as $k$ grows. The most famous illustration of this fact is the excision result due to Goodwillie, Klein and Weiss (\cite{goodwillie_klein},\cite{goodwillie}). Given a smooth $D$-dimensional manifold $W$ with $d+3\leq D$, they prove that the Taylor tower of the presheaf $\Emb(-,W):\cMan_d^{op}\rightarrow\cS$ converges for closed manifolds $M$, that is there is a weak homotopy equivalence

\[
\Emb(M,W)\xrightarrow{\simeq}T_\infty(\Emb(-,W))(M).
\]

This motivates to study more thoroughly the $T_\infty$ term of the tower. To do so, let $O_\infty(M)$ be the poset of open subsets $U\subset M$ diffeomorphic to a finite disjoint union of euclidean balls. By the proof of \cite[Lemma 6.4]{knudsen} (see also \cite[Proposition 8.3]{manifold}), we have an equivalence of disk-presheaves

\[
\colim{U\in O_\infty(M)}{\Emb(-,U)}\simeq\Emb(-,M).
\]

For a presheaf $F:\cMan^{op}_d\rightarrow\cS$, plugging into the definition of $T_\infty F$ provides weak homotopy equivalences 
\begin{align*}
& T_\infty F(M)=\Map_{\Psh(\cDisk_{<\infty})}(\Emb(-,M),F) \\
& \simeq \limi{U\in O_\infty(M)}{\Map_{\Psh(\cDisk_{<\infty})}(\Emb(-,U),F)} \\
& \simeq\limi{U\in O_\infty(M)}{F(U)}
\end{align*}
where the second one comes from the Yoneda lemma. For example, taking $F=\Emb(-,W)$, this recovers the original statement of the excision result as a (homotopy) limit

\[
\Emb(M,W)\simeq\limi{U\in O_\infty(M)}{\Emb(U,W)}.
\]

Albeit this result is of great theoretical interest, in practice the limit is hard to compute since the poset $O(M)_{\infty}$ is not easily tractable. This motivates to find more manageable subposets $\cU\subset O(M)_{\infty}$ for which the limit stays invariant.

One of the first results exploiting this idea appears in \cite{sinha}, where Sinha develops a cosimplicial space to compute the homotopy type of the space of longs knots modulo immersions $\overline{\Emb}^{\partial}(I,I^D)$, with $I=[0,1]$ and $D\geq 4$. His idea can be paraphrased as follows: let $\mathcal{I}$ be the poset of open subsets $U\subset I$ whose complements are of the form 
\[
(I\setminus U)= [a_0,b_0]\cup\cdots\cup[a_p,b_p]
\]
with $0<a_0<b_0<a_1<\cdots<b_{p-1}<a_p<b_p$. On one hand, this poset gives the weak homotopy equivalence 
\[
\overline{\Emb}^{\partial}(I,I^D)\simeq\holim{U\in\mathcal{I}^{op}}{\overline{\Emb}^{\partial}(U,I^D)}.
\] 
On the other hand, the Dwyer-Kan localization of $\mathcal{I}$ at isotopy equivalences provides an equivalence of categories 
\[
\mathcal{I}[\text{isotopies}^{-1}]\xrightarrow{\simeq}\Delta^{op}.
\]
Combining these two results provides a cosimplicial space $\mathcal{K}^{\bullet}$ with 
\[
\overline{\Emb}^{\partial}(I,I^D)\simeq\holim{\Delta}{\mathcal{K}^\bullet}
\] 
and objectwise homotopy equivalences $\mathcal{K}^p\simeq \Conf(p,\R^D)$. This model allows for example to compute the rational homotopy type of $\overline{\Emb}^{\partial}(I,I^D)$, underlining its relevance.

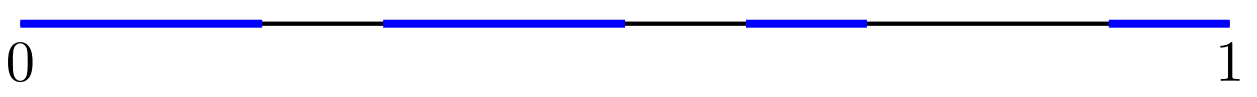
\begin{figure}[ht]
\centering
\begin{tikzpicture}[scale=8]
  \draw[thick] (0,0) -- (1,0);

  \node[below] at (0,0) {$0$};
  \node[below] at (1,0) {$1$};

  \draw[blue, line width=1.5pt] (0,0) -- (0.2,0);
  \draw[blue, line width=1.5pt] (0.3,0) -- (0.5,0);
  \draw[blue, line width=1.5pt] (0.6,0) -- (0.7,0);
  \draw[blue, line width=1.5pt] (0.9,0) -- (1,0);
\end{tikzpicture}
\caption{An open subset $U\in\mathcal{I}$.}
\label{fig:subintervals}
\end{figure}

Expanding on this idea, Boavida de Brito, Lambrechts, Pryor and Songhafouo Tsopméné develop in \cite{model} a cosimplicial model for spaces of embeddings $\Emb(M,W)$ out of a simplicial set $X$ whose geometric realization $\norm{X}$ is homeomorphic to $M$. To do so, they construct a simplicial presheaf  $\Psi_M:\Delta^{op}\rightarrow\Psh(\cDisk_{<\infty})$ such that

\begin{enumerate}
\item $\colim{\Delta^{op}}{\Psi_M}\simeq\Emb(-,M)$,
\item $(\Psi_M)_p\simeq\Emb(-,X_p\times\R^d)$ for all $p\geq 0$.
\end{enumerate}

In particular, each presheaf $(\Psi_M)_p$ is representable. Since the Yoneda embedding
\[
\mathcal{Y}:\cDisk_{<\infty}\hookrightarrow\Psh(\cDisk_{<\infty})
\]
is fully faithful, this implies that the simplicial presheaf $\Psi_M$ factors through $\cDisk_{<\infty}$. In this document, we prove that the simplicial presheaf does not only factor through $\cDisk_{<\infty}$ but also lifts to the $\infty$-category $\cDisk_{/M}$ of disks \emph{embdedded in $M$}.

More precisely, we prove the following statement.

\ \

\begin{theorem}\label{main_simplicial_model}
Let $M$ be a closed $d$-dimensional smooth manifold and $X$ be a finite simplicial set whose geometric realization is homeomorphic to $M$. Then there exists a cofinal functor $\Gamma_X:\Delta^{op}\rightarrow\cDisk_{/M}$ making the following triangle commute:
\[
\begin{tikzcd}
\Delta^{op}\arrow[dr,"X"']\arrow[r,"\Gamma_X"] & \cDisk_{/M}\arrow[d,"\pi_0"] \\
& \Fin
\end{tikzcd}.
\]
\end{theorem}

Albeit the proof of the theorem is done in the case without additional tangential structure,
we will deduce immediately a version for structured manifolds (\Cref{main_structured}).

While the proof relies essentially on the same idea of \cite{model}, its implementation uses different methods. In particular, the results presented in this paper are independant of the ones of \cite{model}. Before giving an overview, we discuss how to recover the simplicial presheaf from \cite{model} and some applications and perspectives suggested in \cite{model}. 

Consider again the Yoneda embedding $\mathcal{Y}:\cDisk_{<\infty}\hookrightarrow\Psh(\cDisk_{<\infty}),U\mapsto\Emb(-,U)$. A formal argument gives an equivalence of disk-presheaves
\[
\colim{U\in\cDisk/M}{\mathcal{Y}(U)}\simeq\Emb(-,M).
\] 

One then gets a simplicial presheaf $\mathcal{Y}\circ\Gamma_X:\Delta^{op}\rightarrow\Psh(\cDisk)$.
Similarly to the one constructed in \cite{model}, this simplicial presheaf satisfies

\begin{enumerate}
\item $\colim{\Delta^{op}}{(\mathcal{Y}\circ\Gamma_X)}\simeq\Emb(-,M)$,
\item $(\mathcal{Y}\circ\Gamma_X)_p\simeq\Emb(-,X_p\times\R^d)$ for all $p\geq 0$
\end{enumerate}

where the first equivalence is a consequence of the cofinality of $\Gamma_X$. Moreover, for a presheaf $F:\cMan^{op}\rightarrow\cS$, one then gets the cosimplicial space
\[
\Omega_F:=\Delta\xrightarrow{\Gamma_X^{op}}\cDisk_{/M}^{op}\rightarrow\cDisk^{op}\hookrightarrow\cMan^{op}\rightarrow\cS.
\]

With this cosimplicial space $\Omega_F$, we can recover \cite[Theorem 1.1]{model}.

\begin{proposition}\label{cosimplicial_space_from_Gamma}
Given a presheaf $F:(\cMan^G)^{op}\rightarrow\cS$, the cosimplicial space $\Omega_F$ satisfies
\begin{enumerate}
\item $T_\infty F(M)\simeq \limi{\Delta}{\Omega_F}$.
\item For each cosimplicial map $\theta:[p]\rightarrow [q]$, there is an embedding $\iota_\theta:X_q\times\R^d\hookrightarrow X_p\times\R^d$ such that $\pi_0(\iota_\theta)\cong (X_\theta:X_q\rightarrow X_p)$ and a weak equivalence of morphisms:
\[(\Omega_F^\theta:\Omega_F^p\rightarrow\Omega_F^q)\simeq (\iota_\theta^*:F(X_p\times\R^d)\rightarrow F(X_q\times\R^d)). 
\] 
\end{enumerate}
\end{proposition}

For a smooth manifold $W$ of dimension at least $d+3$, combining this weak homotopy equivalence and the Goodwillie-Klein-Weiss excision estimate yields a cosimplicial space $\Omega_W$ with 
\[
\Emb(M,W)\simeq\limi{\Delta}{\Omega_W}.
\] 

Following the results of \cite{sinha,sinha_bis,moriya}, one can study the rational cohomology $H^*(\Emb(M,W),\Q)$ using the rational Bousfield-Kan spectral sequence of this cosimplicial space. This is done in \cite{bkss_haya} for spaces of embeddings (modulo immersions) of a parallelizable manifold into an Euclidean space.

Another application concerns factorization homology. Namely, given a symmetric monoidal $\infty$-category $\cC^{\otimes}$ and a $d$-disk algebra $A:\cDisk_{<\infty}\rightarrow\cC^{\otimes}$, the factorization homology of $A$ along $M$ is computed as the colimit 
\[
\int_M A\simeq\colim{\cDisk_{/M}}{A}.
\] 
\Cref{main_simplicial_model} implies that it can be reduced to computing 
\[
\colim{\Delta^{op}}{(A\circ\Gamma_X)}.
\]

Finally, it would also be interesting to extend \Cref{main_simplicial_model} to manifolds with boundaries. 
For a manifold $M$ with nonempty boundary and a monomorphism of simplicial sets $A\hookrightarrow X$ whose geometric realization is homeomorphic to the inclusion $\partial M\hookrightarrow M$, one expects to have a similar model $\Gamma_{X,A}:\Delta^{op}\rightarrow\cDisk^{\partial}_{/M}$, where the superscript indicates a fixed boundary condition between the embeddings.
Arguably, the boundary introduces additional technicalities but the spirit of the proof should stay identical.

\subsection{Overview of the proof}\label{overview_simplicial_model}

The first step of the proof is to replace the $\infty$-category $\cDisk_{/M}$ by the \emph{configuration category} of $M$, $\con(M)$. Following the particle model described in \cite[Section 3.1]{confcat}, this $\infty$-category informally consists of configurations of points in $M$ and paths between them where points can merge but not separate.

The formalization of this idea relies on \emph{exit path categories} as defined in \cite[Appendix A.5]{ha}, which we recall in 
\Cref{section_exit_paths}.  Namely, for each $n\geq 0$, consider $\prt{\us{n}}$ to be the poset of equivalence relations of $\us{n}$ with the partial order given by reverse inclusions. This poset provides a stratification by "diagonals" on $M^n$ where the top stratum is the configuration space $\Conf(n,M)$ and lower strata are defined by whether or not some components agree. The dual of the induced exit path category, $\enter{\prt{\us{n}}}{M^n}$, will then consist of $n$-tuples of points in $M$ and reverse exit paths, i.e paths where points can merge but not separate.

To account for different numbers of points, we then observe that the assignment $\us{n}\mapsto\enter{\prt{\us{n}}}{M^n}$ extends to a contravariant functor $\Fin^{op}\rightarrow\Cat$ from the category of finite sets. By the straightening-unstraightening theorem, this provides a cartesian fibration $\X_M\rightarrow\Fin$. Informally, $\X_M$ is an $\infty$-category whose objects are pairs $(\us{n}\in\Fin,x\in M^n)$ and morphisms $(\us{n},x)\rightarrow (\us{m},y)$ are triangles

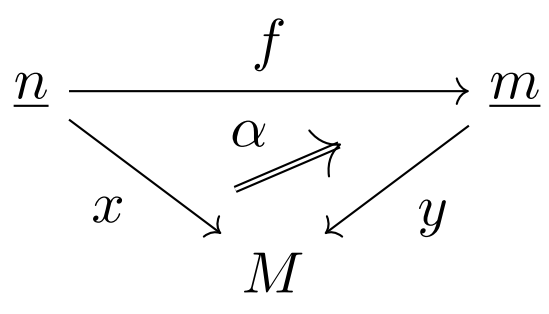
\[
\begin{tikzpicture}[baseline=(B.base)]
\node (A) at (-1.6,1.2) {$\us{n}$};
\node (C) at ( 1.6,1.2) {$\us{m}$};
\node (B) at ( 0,0) {$M$};

\draw[->] (A) -- node[above] {$f$} (C);
\draw[->] (A) -- node[below left] {$x$} (B);
\draw[->] (C) -- node[below right] {$y$} (B);

\draw[double,->] (-0.25,0.55) -- node[above left] {$\alpha$} (0.45,0.85);
\end{tikzpicture}
\]

where $\alpha$ is a reverse exit path from $x$ to $yf$ in $M^n$. The configuration category $\con(M)$ is then defined as the full subcategory whose objects are pairs $(\us{n},x)$ with $x\in \Conf(n,M))$. In \Cref{section_conf_cat}, we prove that $\con(M)$ is equivalent to a localization of $\X_M$ and prove the equivalence of $\infty$-categories $\con(M)\simeq\cDisk_{/M}$ by adapting a result of \cite[Section 3.2]{confcat}.

\begin{figure}[ht]
\centering
\begin{tikzpicture}[scale=1.2, line cap=round, line join=round]
  \usetikzlibrary{decorations.markings}

  \def\A{3}      
  \def\B{1.4}    
  \def\C{1.2}    
  \def\D{0.5}    

  \tikzset{
    midarrow/.style={
      postaction={decorate},
      decoration={
        markings,
        mark=at position 0.5 with {\arrow{>}}
      }
    }
  }

  \fill[gray!20, even odd rule]
    (0,0) ellipse ({\A} and {\B})
    (0,0) ellipse ({\C} and {\D});

  \draw[dashed] (\A,0)
    arc[start angle=0,end angle=180,x radius=\A,y radius=\B];
  \draw[thick] (-\A,0)
    arc[start angle=180,end angle=360,x radius=\A,y radius=\B];

  \draw[dashed] (\C,0)
    arc[start angle=0,end angle=180,x radius=\C,y radius=\D];
  \draw[thick] (-\C,0)
    arc[start angle=180,end angle=360,x radius=\C,y radius=\D];


  \coordinate (B1) at (-2.05,0.52);   
  \coordinate (B2) at (-1.15,1.02);
  \coordinate (B3) at (-1.85,-0.78);  

  \coordinate (R1) at ( 2.00,0.56);   
  \coordinate (R2) at ( 2.10,-0.46);

  \fill[blue] (B1) circle (1.2pt);
  \fill[blue] (B2) circle (1.2pt);
  \fill[blue] (B3) circle (1.2pt);

  \node[blue, above left] at (B1) {\scriptsize $1$};
  \node[blue, left] at (B2) {\scriptsize $2$};
  \node[blue, above left] at (B3) {\scriptsize $3$};

  \fill[red] (R1) circle (1.2pt);
  \fill[red] (R2) circle (1.2pt);

  \node[red, above right] at (R1) {\scriptsize $1$};
  \node[red, below right] at (R2) {\scriptsize $2$};

  \draw[thick, midarrow] (B1) to[out=18,in=162] (R1);
  \draw[thick, midarrow] (B2) to[out=8,in=130] (R1);
  \draw[thick, midarrow] (B3) to[out=-5,in=-165] (R2);

\end{tikzpicture}
\caption{Depicting a morphism in the configuration category of the torus.}
\end{figure}

\ \

We can make an analogous construction for the interval $I=[0,1]$, taking the linear order into account. 
More precisely, $\con^{<}(I)$ is an $\infty$-category whose objects are configurations $(0<t_1<\cdots<t_p<1)$ in the interval $I=[0,1]$ and morphisms are reverse exit paths. The superscript indicates that configurations are given by linearly ordered points $(0<t_1<\cdots<t_p<1)$ in $I$. We replace the simplex category $\Delta^{op}$ by the \emph{category of finite linear intervals} $\nabla\cong\Delta^{op}$ (\cite{wraith}). Then, the configurations of $\con^<(I)$ are obtained as linearly ordered maps $u:\la p\ra\rightarrow I$ for some object $\la p\ra\in\nabla$. This gives a projection functor $\con^<(I)\rightarrow\nabla\simeq\Delta^{op}$ and we prove it is an equivalence of $\infty$-categories (\Cref{O_equiv_nabla}). 

The construction of the $\infty$-category $\con^<(I)$ is similar to the one of $\con(M)$ following the table of analogies below.

\ \

\begin{table}[ht]
\centering

\renewcommand{\arraystretch}{1.35}

\begin{tabular}{|c|c|}
\hline

$\nabla\cong\Delta^{op}$
&
$\Fin$
\\
\hline

$\Delta^p\cong\{u:\la p\ra\rightarrow I\}$.
&
$M^n\cong\{x:\us{n}\rightarrow M\}$.
\\
\hline

\shortstack{
$\pos{p}$ poset of faces\\
$\la p\ra\twoheadrightarrow\la q\ra$.
}
&
\shortstack{
$\prt{\us{n}}$ poset of equivalence\\
relations on $\us{n}$.
}
\\
\hline

$\enter{\pos{-}}{\Delta^-}:\nabla^{op}\rightarrow\Cat$.
&
$\enter{\prt{-}}{M^-}:\Fin^{op}\rightarrow\Cat$.
\\
\hline

$\F\rightarrow\nabla$ the unstraightening.
&
$\X_M\rightarrow\Fin$ the unstraightening.
\\
\hline

\shortstack{
$\con^<(I)\subset\F$\\
subcategory of configurations.
}
&
\shortstack{
$\con(M)\subset\X_M$\\
subcategory of configurations.
}
\\
\hline

$\con^<(I)\simeq\nabla$.
&
$\con(M)\simeq\cDisk_{/M}$.
\\

\hline
\end{tabular}
\end{table}

\ \

\begin{example}
The following picture shows a morphism from the configuration $u=(0<0.2<0.6<1)$ in blue to the configuration $v=(0<0.5<1)$ in red.
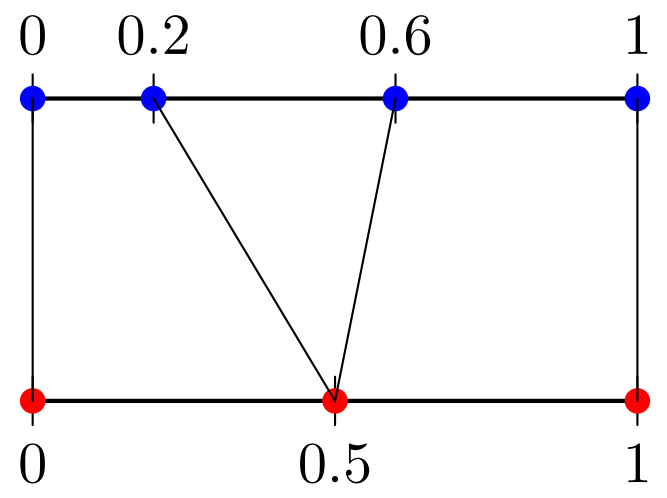
\begin{center}
\begin{tikzpicture}[scale=4, line cap=round, line join=round]

\draw[thick] (0,0) -- (1,0);
\foreach \x/\lab in {0/0, 0.2/0.2, 0.6/0.6, 1/1}{
  \draw (\x,0.04) -- (\x,-0.04);     
  \fill[blue] (\x,0) circle (0.6pt);       
  \node[above] at (\x,0.04) {$\lab$};
}

\def\dy{-0.5}
\draw[thick] (0,\dy) -- (1,\dy);
\foreach \x/\lab in {0/0, 0.5/0.5, 1/1}{
  \draw (\x,\dy+0.04) -- (\x,\dy-0.04);  
  \fill[red] (\x,\dy) circle (0.6pt);         
  \node[below] at (\x,\dy-0.04) {$\lab$};
}

\draw (0,0) -- (0,\dy);
\draw (0.2,0) -- (0.5,\dy);
\draw (0.6,0) -- (0.5,\dy);
\draw (1,0) -- (1,\dy);

\end{tikzpicture}
\end{center}
\end{example}

Because of the equivalences $\con^<(I)\simeq\Delta^{op},\con(M)\simeq\cDisk_{/M}$, to define the relevant functor $\Gamma_X:\Delta^{op}\rightarrow\cDisk_{/M}$, we instead define a functor $\gamma:\con^<(I)\rightarrow\con(M)$. This is the content of \Cref{section_gamma}. It is enough to define a functor $\bar{\gamma}:\F\rightarrow\X_M$ and check it restricts to a functor $\con^<(I)\rightarrow \con(M)$. The construction follows the idea of \cite{model}, replacing subintervals of $I$ by configurations in $I$.

Namely, given an object $u=(0\leq t_1\leq \cdots\leq t_p\leq 1)\in\F$ and an integer $q\geq 0$, we get a finite subset $u^q\subset I^q$ which in turn restricts to a finite subset $\eta_q(u)\subset\Delta^q$ (see the pictures below). Taking the colimit of the geometric simplices along the simplex category $\Delta(X)$ then gives a subset $\eta(u)\subset\norm{X}$. When $u=(0< t_1<\cdots<t_p<1)$ is an element of $\con^<(I)$, we prove that this subset $\eta(u)$ has cardinality $\# X_p$, that is, is an element of $\Conf(X_p,\norm{X})$.

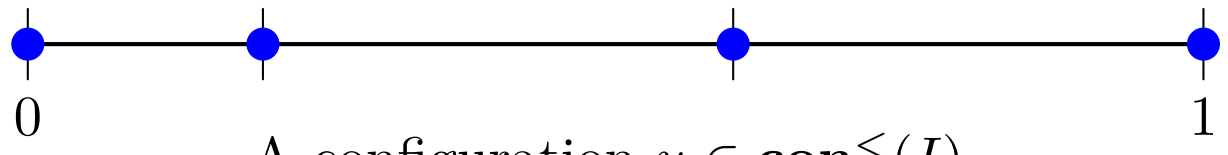
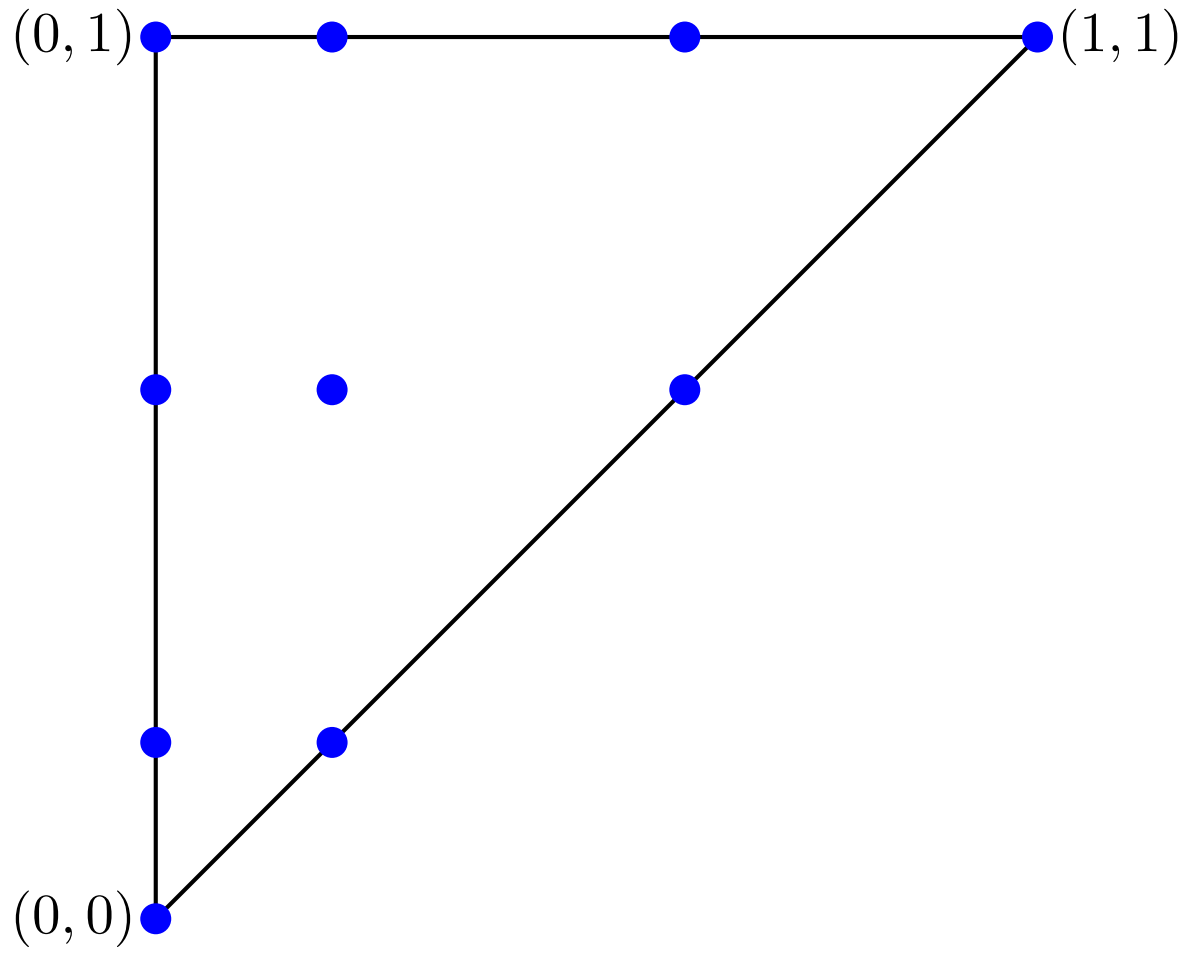
\begin{figure}[H]
\centering
\resizebox{!}{0.70\textheight}{%
\begin{tikzpicture}

\tikzset{
  onearrow/.style={
    postaction={decorate},
    decoration={
      markings,
      mark=at position 0.5 with {\arrow{>}}
    }
  },
  twoarrows/.style={
    postaction={decorate},
    decoration={
      markings,
      mark=at position 0.46 with {\arrow{>}},
      mark=at position 0.54 with {\arrow{>}}
    }
  }
}

\begin{scope}[yshift=0cm]
\begin{scope}[scale=8, line cap=round, line join=round]
  \draw[thick] (0,0) -- (1,0);

  \foreach \x/\lab in {0/0, 0.2/0.2, 0.6/0.6, 1/1} {
    \draw (\x,0.03) -- (\x,-0.03);
    \fill[blue] (\x,0) circle (0.4pt);
  }
  
  \node[below] at (0,-0.03) {$0$};
  \node[below] at (1,-0.03) {$1$};
\end{scope}

\node[align=center] at (4,-0.8)
  {A configuration $u\in\con^<(I)$.};
\end{scope}

\begin{scope}[xshift=1cm,yshift=-8cm]
\begin{scope}[scale=6, line cap=round, line join=round]

  \draw[thick] (0,0) -- (0,1) -- (1,1) -- cycle;

  \foreach \x in {0,0.2,0.6,1} {
    \foreach \y in {0,0.2,0.6,1} {
      \ifdim\y pt<\x pt
      \else
        \fill[blue] (\x,\y) circle (0.5pt);
      \fi
    }
  }

  \node[left]  at (0,0) {$(0,0)$};
  \node[left]  at (0,1) {$(0,1)$};
  \node[right] at (1,1) {$(1,1)$};

\end{scope}

\node[align=center] at (3,-0.8)
  {The configuration $\eta_2(u)$ in $\Delta^2$.};
\end{scope}

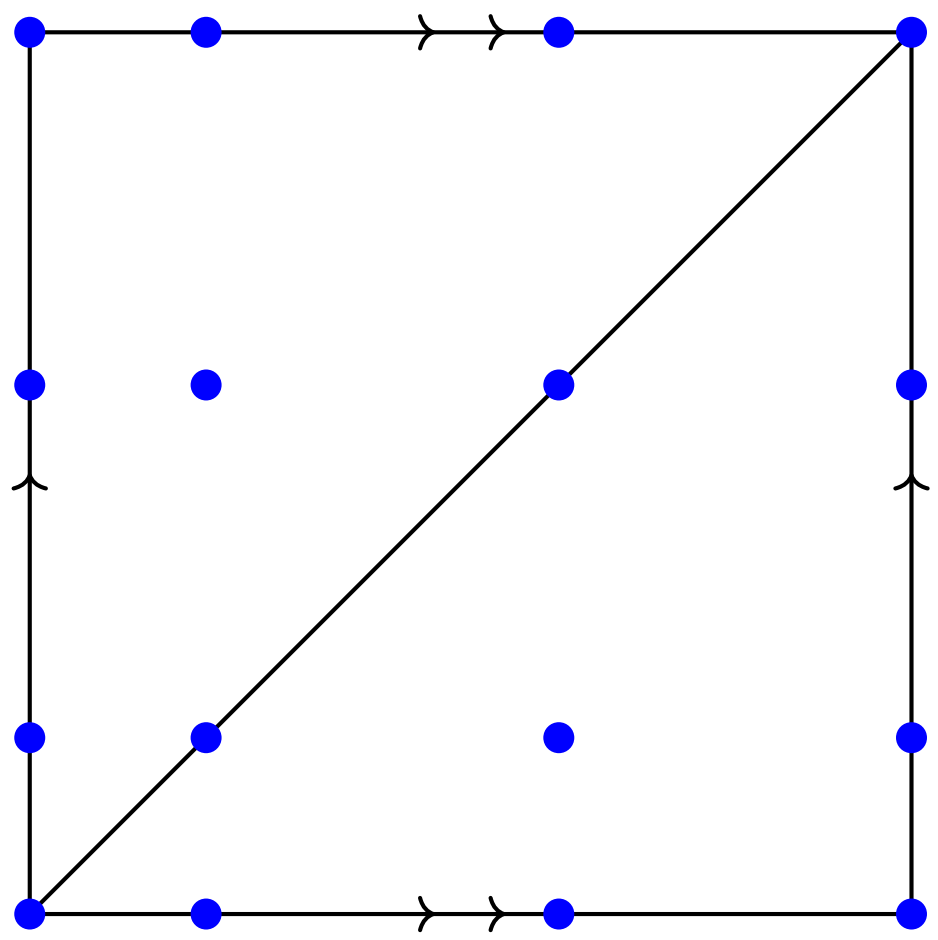
\begin{scope}[xshift=1cm,yshift=-15.8cm]
\begin{scope}[scale=6, line cap=round, line join=round]

  \draw[thick, twoarrows] (0,0) -- (1,0); 
  \draw[thick, onearrow]  (0,0) -- (0,1); 
  \draw[thick, onearrow]  (1,0) -- (1,1); 
  \draw[thick, twoarrows] (0,1) -- (1,1); 

  \draw[thick] (0,0) -- (1,1);

  \foreach \x in {0,0.2,0.6,1} {
    \foreach \y in {0,0.2,0.6,1} {
      \fill[blue] (\x,\y) circle (0.5pt);
    }
  }

\end{scope}

\node[align=center, text width=11cm] at (3,-0.9)
  {The configuration $\eta(u)$ in $\norm{X}$ for
  $X=(\Delta[1]/\partial\Delta[1])^2$.};
\end{scope}

\end{tikzpicture}
}
\end{figure}

This construction sends reverse exit paths in $\F$ to reverse exit paths in $\X_M$, thereby providing the desired functor $\bar{\gamma}:\F\rightarrow\X_M$.

In fact, the construction is compatible with a finer stratification on the topological spaces $\{\norm{X^n}\}_{n\geq 0}$ than the one by diagonals, that is reviewed in \Cref{section_stratification_realization}.
The geometric realization $\norm{X}$ admits a stratification by its simplices. Namely, the poset $\pos{X}$ consisting of nondegenerate simplices of $X$ and face relations between them provides a stratification $\norm{X}\rightarrow\pos{X}$. This construction extends to a functor 
\[
\Fin^{op}\rightarrow\Cat,\us{n}\mapsto\enter{\pos{X^n}}{\norm{X^n}}
\]
which corresponds to a cartesian fibration $\B_X\rightarrow\Fin$ under the straightening-unstraightening equivalence. It turns out that for any $\us{n}\in\Fin$, this $\pos{X^n}$-stratification is a refinement of the $\prt{\us{n}}$-stratification on $\norm{X^n}\cong M^n$. This translates into a cartesian functor $\B_X\rightarrow\X_M$ and the functor $\bar{\gamma}:\F\rightarrow\X_M$ factors through this refinement via a functor $\delta:\F\rightarrow\B_X$.

This refinement appears to be useful when proving the cofinality of the functor $\gamma:\con^<(I)\rightarrow \con(M)$. This is the content of \Cref{section_cofinality}. Combining the fact that the inclusions $\con^<(I)\hookrightarrow\F,\con(M)\hookrightarrow\X_M$ are fully faithful right adjoints and a localization result about exit path categories, we reduce the proof of the cofinality for $\gamma$ to proving that the functor $\delta:\F\rightarrow\B_X$ is cofinal. 

By Quillen's Theorem A, it is enough to show that for any object $(x:\us{n}\rightarrow M)\in\B_X$, the slice category $\delta_{x/}$ is weakly contractible. Informally, the $\infty$-category $\delta_{x/}$ has objects reverse exit paths $(x\rightarrow\delta(u))$ in $\enter{\pos{X^n}}{\norm{X^n}}$, with $u\in\F$, and morphisms being homotopies between reverse exit paths in $\enter{\pos{X^n}}{\norm{X^n}}$

\[
\begin{tikzpicture}[baseline=(B.base)]
\node (A) at (-1.6,0) {$\delta(u)$};
\node (C) at ( 1.6,0) {$\delta(v)$};
\node (B) at ( 0,1.2) {$x$};

\draw[->] (A) -- node[below] {$\delta(\alpha)$} (C);
\draw[->] (B) -- node[below left] {} (A);
\draw[->] (B) -- node[below right] {} (C);

\draw[double,->] (-0.25,0.55) -- node[above left] {} (0.35,0.70);
\end{tikzpicture}
\]

where $\alpha:u\rightarrow v$ is a morphism in $\F$.

\Cref{balancing_lemma} implies that it amounts to prove the cofinality of functors 
\[
\lambda_n:\cE_n\rightarrow\enter{\pos{X^n}}{\norm{X^n}},n\geq 0.
\]

The $\infty$-category $\cE_n$ has objects pairs $(f:\us{n}\rightarrow X_p,u:\la p\ra\hookrightarrow I)$ for some $p\geq 0$ and $\lambda_n$ sends such an object to the composition
\[
\us{n}\xrightarrow{f}X_p\xrightarrow{\eta(u)}\norm{X}.
\]

This $\infty$-category comes with a cartesian fibration $\cE_n\rightarrow\Delta(X^n)^{op}$ over the dual of the category $\Delta(X^n)^{op}$ of simplices of $X^n$. Its straightening is the functor $\enter{\pos{-}}{\Delta^-}:\Delta(X^n)\rightarrow\Cat$. \Cref{stratified_realization} identifies $\enter{\pos{X^n}}{\norm{X^n}}$ as the colimit of this diagram, so that $\lambda_n$ is the localization functor exhibiting this colimit, hence is cofinal.

The following diagram summarizes the different $\infty$-categories and functors involved in the construction.
\[
\begin{tikzcd}
\con^<(I)\arrow[d,hook]\arrow[rr,"\gamma"] & & \con(M)\arrow[d,hook] \\
\F\arrow[d]\arrow[r,"\delta"] & \B\arrow[r,"\mu"]\arrow[dr] & \X_M\arrow[d] \\
\nabla\arrow[rr,"X"'] & & \Fin  
\end{tikzcd}.
\]

\subsection{Conventions and notations}
We establish here our conventions for classical objects from the literature. Objects specific to the document will be defined in the relevant sections.

\begin{itemize}
\item We write $\Delta$ for the classical simplicial category whose objects are linearly ordered sets $[p]=\{0\leq 1\leq\cdots\leq p\}$ for $p\geq 0$.
\item $\Fin$ is the ordinary category whose objects are finite sets $\us{n}:=\{1,\cdots,n\}$ for each $n\geq 0$ and morphisms are maps between them.
\item $\Top$ is the ordinary category of topological spaces and $\sSet$ the one of simplicial sets.
\item  $\Delta^p=\{(t_1,\cdots,t_p)\in [0,1]^p:0\leq t_1\leq\cdots\leq t_p\leq 1\}$ is the geometric simplex and $\Delta[p]=\Delta(-,[p])$ is the representable simplicial set.
\item Given a simplicial set $X$, we write $\Delta(X)$ for its category of elements. For a simplex $\sigma:\Delta[p]\rightarrow X$, we write $\Delta[\sigma]:=\Delta[p],\Delta^\sigma:=\Delta^p$ for the (geometric) simplex it indexes.
\item We use \emph{quasicategories} as a model for $(\infty,1)$-categories, that is simplicial sets satisfying the inner horn-filling condition. We will simply call them $\infty$-categories. We insist on the fact that our construction is heavily dependent on this specific model.
\item We consider all ordinary categories to be locally small. To distinguish them from $\infty$-categories, we will keep calling them ordinary categories. We will sometimes safely abuse the notation by considering an ordinary category as an $\infty$-category via its nerve. We write $\text{Cat}$ the ordinary category of small ordinary categories. The nerve functor is denoted $N:\text{Cat}\rightarrow\sSet$.
\item Given $\infty$-categories $\cC,\cD$, we write $\Fun(\cC,\cD)$ to be the $\infty$-category given by $\Fun(\cC,\cD)_p=\sSet(\cC\times\Delta[p],\cD)$.
\item We use \emph{Kan complexes} as a model for $\infty$-groupoids. We will call them \emph{spaces}.
\item Given an ordinary category $\cC$ (respectively an $\infty$-category) and objects $x,y\in\cC$, we write $\cC(x,y)$ for the set (respectively space) of morphisms between them.
\item Given an $\infty$-category $\cC$, we write $\cC^{\simeq}$ its maximal subgroupoid and $\norm{\cC}$ its classifying space.
\item We write $\text{Kan}$ for the simplicial category having objects Kan complexes and mapping spaces $\Fun(\cC,\cD)$. We model the $\infty$-category of (small) spaces by $\cS:=N^h(\text{Kan})$ where $N^h$ is the coherent nerve functor.
\item We write $\mathcal{Q}\text{Cat}$ for the simplicial category having objects quasicategories and mapping spaces $\Fun(\cC,\cD)^{\simeq}$. We model the $\infty$-category of (small) $\infty$-categories by $\Cat:=N^h(\mathcal{Q}\text{Cat})$. We write $\Cat^{\delta}$ for the ordinary category having objects $\infty$-categories and functors between them (so composition is strict).
\item We take the convention that "cofinal" means colimit preserving.
\item The hooked arrow $"\hookrightarrow"$ will either designate an injective map, an embedding between manifolds or an inclusion of full subcategory.
\end{itemize}

\subsection*{Acknowledgements}
This work was realized during the course of the author's PhD program.

He thanks his advisors Pascal Lambrechts and Nathalie Wahl for their careful reading and for providing valuable comments on earlier versions of this document.

He also thanks Pascal Lambrechts for generously sharing a preliminary version of \cite{model}, on which the general idea underlying this paper is based.

The author was supported by the "Fonds de la Recherche Scientifique–FNRS" under the Aspirant Fellowship FC 49827. 

\subsection*{Use of generative AI}
During the preparation of this document, generative AI has been used for skimming more rapidly through the relevant literature, for reviewing early versions by pointing out mistakes and inconsistencies and assisting with the document layout. The models used were Gemini 3.1,3.5,3.6 and ChatGPT-5.6.

\section{About $\infty$-categories}

This section is devoted to quick recollections about $\infty$-categories and proving related results. We use the setting of \emph{quasicategories} as initiated by Joyal and developed by Lurie in \cite{htt} (see also \cite{kerodon}). 

\subsection{Cartesian fibrations}
Cartesian fibrations and the unstraightening-straightening theorem are the $\infty$-categorical analog of Grothendieck fibrations and the Grothendieck construction. Namely, given an $\infty$-category $\cC$ and a functor $F:\cC^{op}\rightarrow\Cat$, there exists a cartesian fibration $p:\cE\rightarrow\cC$ parametrizing $F$. Informally speaking, $\cE$ has objects consisting of pairs $(c\in\cC,x\in F(c))$, morphisms $(c,x\in F(c))\rightarrow (d,y\in F(d))$ consisting of pairs $(f:c\rightarrow d\in\cC,\alpha:x\rightarrow F(f)(y)\in F(c))$ and $p$ is the obvious projection. We refer readers unfamiliar with cartesian fibrations to \cite[2.4]{htt}. Conversely, given a cartesian fibration $\cE\rightarrow\cC$, it parametrizes an essentially unique functor $F:\cC^{op}\rightarrow\Cat$. This is the content of the famous \emph{straightening-unstraightening} theorem. In the following, we write $\cart_{/\cC}$ for the $\infty$-category of cartesian fibrations over $\cC$.
 
\begin{theorem}\cite[Theorem 3.2.0.1]{htt}
Let $\cC$ be an $\infty$-category. There is an equivalence of $\infty$-categories $\text{St}_{\cC}:\cart_{/\cC}\simeq\Fun(\cC^{op},\Cat):\text{Un}_{\cC}$ through the \emph{straightening functor} $\text{St}_{\cC}$ and the \emph{unstraightening functor} $\text{Un}_{\cC}$.
\end{theorem}

Whenever $\cC$ is an ordinary category, given a strict functor $F:\cC^{op}\rightarrow\Cat^{\delta}$, its unstraightening admits an explicit description through the \emph{weighted nerve}.

\begin{definition}\cite[\href{https://kerodon.net/tag/025X}{Definition 5.3.3.1}]{kerodon}
Let $\cC$ be an ordinary category and $F:\cC^{op}\rightarrow\Cat^{\delta}$ be a functor between ordinary categories. Its \emph{weighted nerve} is the simplicial set $N^F$ which has $p$-simplices pairs of a chain $(c_0\xleftarrow{f_1} c_1\xleftarrow{f_2}\cdots\xleftarrow{f_p} c_p)\in N(\cC)_p$ and a commutative diagram
\[
\begin{tikzcd}
\Delta[0]\arrow[r,"d^0"]\arrow[d] & \Delta[1]\arrow[r,"d^0"]\arrow[d] & \cdots\arrow[r,"d^0"]\arrow[d] & \Delta[p]\arrow[d] \\
F(c_0)\arrow[r,"F(f_1)"'] & F(c_1)\arrow[r,"F(f_2)"'] & \cdots\arrow[r,"F(f_p)"'] & F(c_p)
\end{tikzcd}.
\]
\end{definition}

\begin{proposition}\cite[\href{https://kerodon.net/tag/046Y}{Corollary 5.3.3.16}]{kerodon}
Let $\cC$ be an ordinary category and $F:\cC^{op}\rightarrow\Cat^{\delta}$ be a functor between ordinary categories. Then its weighted nerve $N^F(\cC)$ is an $\infty$-category and the obvious projection $p^F:N^F(\cC)\rightarrow N(\cC)$ is a cartesian fibration equivalent to the unstraightening of the functor $\cC^{op}\xrightarrow{F}\Cat^{\delta}\rightarrow\Cat$.
\end{proposition}

By abuse of notation, we will sometimes call the cartesian fibration $p^F:N^F(\cC)\rightarrow N(\cC)$ the weighted nerve.

In this document, all unstraightenings will come from strict functors $\cC^{op}\rightarrow\Cat^{\delta}$ hence will be modelled using the weighted nerve. Note that there is also a covariant version of the weighted nerve for strict functors $\cC\rightarrow\Cat^\delta$.

\begin{lemma}\label{core_cartesian_fib}
Let $F:\cC^{op}\rightarrow\Cat$ be a functor with unstraightening the cartesian fibration $p:\cE\rightarrow\cC$. Then the composition $(\cC^{op})^{\simeq}\hookrightarrow\cC^{op}\xrightarrow{F}\Cat\xrightarrow{(-)^{\simeq}}\cS$ has unstraightening the right fibration $\cE^{\simeq}\rightarrow\cC^{\simeq}$.
\end{lemma}

\begin{proof}
Let $R\rightarrow\cC$ be the right fibration being the unstraightening of the composition $\cC^{op}\rightarrow\Cat\xrightarrow{(-)^{\simeq}}\cS$. Because of the naturality of the inclusion of the maximal subgroupoid inside an $\infty$-category, this provides a cartesian functor
\[
\begin{tikzcd}
R\arrow[rr]\arrow[dr] & & \cE\arrow[dl] \\
& \cC &
\end{tikzcd}
\]
whose fiber over an object $c\in\cC$ is the inclusion $F(c)^{\simeq}\hookrightarrow F(c)$.

Now consider the pullback
\[
\begin{tikzcd}
K\arrow[d]\arrow[r]\arrow[dr,very near start,"\lrcorner",phantom] & R\arrow[d] \\
\cC^{\simeq}\arrow[r,hook] & \cC
\end{tikzcd}.
\]
The functor $K\rightarrow\cC^{\simeq}$ is a right fibration over a space hence $K$ is itself a space. Thus the inclusion $K\hookrightarrow R\hookrightarrow\cE$ factors through $\cE^{\simeq}$. Now given an object $c\in\cC$, the inclusion of fiber $F(c)\hookrightarrow\cE$ restricts to an inclusion $F(c)^\simeq\hookrightarrow\cE$ which factors through an inclusion $F(c)^\simeq\hookrightarrow\cE^{\simeq}$. We thus get the following ladder 

\[
\begin{tikzcd}
F(c)^{\simeq}\arrow[d,hook]\arrow[r,equal] & F(c)^\simeq\arrow[d,hook] \\
K\arrow[r]\arrow[d] & \cE^\simeq \arrow[d] \\
\cC^{\simeq}\arrow[r,equal] & \cC^{\simeq}
\end{tikzcd}.
\]

The left side is a homotopy fiber sequence by definition of the right fibration $K\rightarrow\cC^{\simeq}$. Moreover the right side is obtained by applying the core functor to the fiber sequence $F(c)\hookrightarrow\cE\rightarrow\cC$. Because the core functor is a right adjoint, it preserves limits hence the right side is also a homotopy fiber sequence.

Because the ladder is a map between homotopy fiber sequences with weakly homotopic fibers and bases, we deduce that the map $K\rightarrow\cE^{\simeq}$ is a weak homotopy equivalence too. Thus the right fibrations $(K\rightarrow\cC^{\simeq}),(\cE^{\simeq}\rightarrow\cC^{\simeq})$ and therefore their straightenings are equivalent. But the straightening of the right fibration $K\rightarrow\cC^\simeq$ is equivalent to the composition $(\cC^{op})^\simeq\hookrightarrow\cC^{op}\xrightarrow{F}\Cat\xrightarrow{(-)^\simeq}\cS$.
\end{proof}

\subsection{Adjunctions}
Two functors $F:\cC\rightleftarrows\cD:G$ between $\infty$-categories form an adjunction \(F\dashv G\) if
there is an equivalence of mapping spaces \(\cD(F(x),y)\simeq\cC(x,G(y))\), natural in \(x\in\cC,\
  y\in\cD\). In practice adjoints are produced via the
comma-category criterion generalizing the classical one: \(F\)
admits a right adjoint if and only if for every \(y\in\cD\) the slice \(F_{/y}\)  has a final object and then the counit $F(G(y))\rightarrow y$ is
that final object. In practice, we have the two following criteria to detect adjoint functors.

\begin{proposition}\label{local_unit}\cite[Remark 5.1.11]{land}
Let $G:\cD\rightarrow\cC$ be a functor between $\infty$-categories. Assume that for each object $x\in\cC$, there exists an object $Fx\in\cD$ and a morphism $\eta_x:x\rightarrow G(Fx)$ such that the composition
\[
\cD(Fx,y)\rightarrow\cC(G(Fx),G(y))\xrightarrow{\eta_x^*}\cC(x,G(y))
\] 
is a weak homotopy equivalence for all $y\in\cD$. Then the assignment $x\mapsto Fx$ extends to a functor $F:\cC\rightarrow\cD$ which is left adjoint to $G$. Moreover, the assignment $x\mapsto\eta_x$ defines the unit of the adjunction.
\end{proposition}

\begin{proposition}\label{adjoint_slice}\cite[Proposition 5.2.6]{land}
Let $F:\cC\rightarrow\cD$ be a functor between $\infty$-categories. Then $F$ admits a right adjoint if and only if for all $x\in\cD$, the slice $F_{/x}$ admits a terminal object.
\end{proposition}

\begin{lemma}\label{pullback_adj}
Consider an adjunction

\[
\begin{tikzcd}
\cC \arrow[r, shift left=1.5ex, "F"]
& 
\cD \arrow[l, shift left=1.5ex, "G"]
\arrow[phantom, from=1-1, to=1-2, "\dashv" rotate=-90, description]
\end{tikzcd}.
\]

For any object $x\in\cC$, the following square is a pullback:

\[
\begin{tikzcd}
\cD_{F(x)/}\arrow[d]\arrow[r]\arrow[dr, phantom, "\lrcorner", very near start] & \cC_{x/}\arrow[d] \\
\cD\arrow[r,"G"'] & \cC
\end{tikzcd}.
\]
\end{lemma}

\begin{proof}
It suffices to prove that there is an equivalence of categories $\mathcal{D}_{F(x)/}\simeq G_{x/}$. Because both functors $\mathcal{D}_{F(x)/}\rightarrow\mathcal{D}, G_{x/}\rightarrow\cD$ are left fibrations, it suffices to prove that their unstraightenings are equivalent.

The unstraightening of the first one is given by the functor $\cD(F(x),-):\cD\rightarrow\cS$ while the second one is $\cC(x,G(-)):\cD\rightarrow\cS$. Because of the adjunction, these are indeed equivalent.
\end{proof}

\begin{lemma}\label{adjoint_cart}
Let $q:\mathcal{G}\rightarrow I$ be a cartesian fibration and assume that $I$ has an initial object $i$. Then the inclusion $j:q^{-1}(i)\hookrightarrow\mathcal{G}$ admits a right adjoint.
\end{lemma}

\begin{proof}
By \cite[Proposition 5.2.6]{land}, it suffices to prove that for any object $y\in\mathcal{G}$, the pullback $j_{/y}$ has a terminal object. Let $f:i\rightarrow q(y)$ be a morphism in $I$ and $\hat{f}:f^*y\rightarrow y$ a cartesian lift of $f$. We claim that $\hat{f}$ is a terminal object of $j_{/y}$.

Let $g:a\rightarrow y$ an object of $j_{/y}$. In particular we have $q(a)=i$. By definition of $j_{/y}$, there is a pullback of mapping spaces

\[
\begin{tikzcd}
j_{/y}(g,\hat{f})\arrow[d]\arrow[r,"j"] & \mathcal{G}_{/y}(g,\hat{f})\arrow[d] \\
q^{-1}(i)(a,f^*y)\arrow[r,"j"'] & \mathcal{G}(a,f^*y)
\end{tikzcd}.
\] 
On one hand, because $q:\mathcal{G}\rightarrow I$ is a cartesian fibration, we have a homotopy fiber sequence

\[
q^{-1}(i)(a,f^*y)\xrightarrow{j}\mathcal{G}(a,f^*y)\rightarrow I(q(a),q(f^*y))
\]

over the identity in $I(i,i)$. Because $i$ is an initial object, we have $I(q(a),q(f^*y))\simeq\ast$ and thus we get a weak homotopy equivalence $ q^{-1}(i)(a,f^*y)\xrightarrow{\simeq}\mathcal{G}(a,f^*y)$.

To conclude the proof, it remains to show that $\mathcal{G}_{/y}(g,\hat{f})$ is weakly contractible.

Because $\hat{f}$ is a cartesian morphism, we have a pullback

\[
\begin{tikzcd}
\mathcal{G}(a,f^*y)\arrow[d]\arrow[r,"\hat{f}"] & \mathcal{G}(a,y)\arrow[d] \\
I(i,i)\arrow[r,"f"'] & I(i,q(y))
\end{tikzcd}.
\]

Both lower mapping spaces are weakly contractible because $i$ is an initial object. Thus we get a weak homotopy equivalence $\mathcal{G}(a,f^*y)\xrightarrow{\hat{f}}\mathcal{G}(a,y)$. In particular the fiber over $g\in\mathcal{G}(a,y)$ is nonempty so that we can pick an element $\alpha:a\rightarrow f^*y$ in it.

Because $\hat{f}$ is a cartesian morphism, we have another pullback

\[
\begin{tikzcd}
\mathcal{G}_{/f^*y}\arrow[r,"\hat{f}"]\arrow[d] & \mathcal{G}_{/y}\arrow[d] \\
I_{/i}\arrow[r,"f"'] & I_{/q(y)}
\end{tikzcd}.
\]

Which gives a pullback on mapping spaces

\[
\begin{tikzcd}
\mathcal{G}_{/f^*y}(\alpha,1_{f^*y})\arrow[d]\arrow[r,"\hat{f}"] & \mathcal{G}_{/y}(g,\hat{f})\arrow[d] \\
I_{/i}(1_i,1_i)\arrow[r,"f"'] & I_{/q(y)}(f,f).
\end{tikzcd}
\]

Once again, because $i$ is an initial object, the lower mapping spaces are weakly contractible. This gives a weak homotopy equivalence $\mathcal{G}_{/f^*y}(\alpha,1_{f^*y})\xrightarrow{\hat{f}} \mathcal{G}_{/y}(g,\hat{f})$. But the left mapping space is weakly contractible because $1_{f^*y}$ is terminal in $\mathcal{G}_{/f^*y}$. Therefore $\mathcal{G}_{/y}(g,\hat{f})$ is also weakly contractible.

\end{proof}

\begin{rem}
Let $F:I^{op}\rightarrow\Cat$ be a functor where $I$ is a small $\infty$-category, not necessarily with an initial object. Consider its colimit cocone $F^\rhd:(I^\lhd)^{op}\simeq (I^{op})^\rhd\rightarrow\Cat$ and let $\bar{q}:\bar{\mathcal{G}}\rightarrow I^\lhd$ be its unstraightening. Since the cone $I^\lhd$ has an initial object $-\infty$, \Cref{adjoint_cart} asserts that the inclusion $\colim{I^{op}}{F}\simeq \bar{q}^{-1}(-\infty)\hookrightarrow\bar{\mathcal{G}}$ has a right adjoint $\rho:\bar{\mathcal{G}}\rightarrow\bar{q}^{-1}(-\infty)$. Now let $q:\mathcal{G}\rightarrow I$ be the unstraightening of $F$. Then, \cite[Theorem 7.4.5.13]{kerodon} implies that the composition $\mathcal{G}\hookrightarrow\bar{\mathcal{G}}\xrightarrow{\rho}\bar{q}^{-1}(-\infty)$ is the localization functor exhibiting $\bar{q}^{-1}(-\infty)\simeq\colim{I^{op}}{F}$ as the localization of $\mathcal{G}$ at cartesian morphisms.
\end{rem}

\subsection{Cofinality}

Let $F:\cC\rightarrow\cD$ be a functor between $\infty$-categories. If $F$ is \emph{cofinal}, $\cD$-shaped colimits can be computed by reduction via $F$ to $\cC$-shaped colimits. In practice, there are various results to check that a functor is cofinal. For example, right adjoints and localizations (which we discuss in the next section) are cofinal. In addition, there is the famous Quillen's Theorem A.

\begin{theorem}\cite[Theorem 4.1.3.1]{htt}
Let $F:\cC\rightarrow\cD$ be a functor between $\infty$-categories. Then $F$ is cofinal if and only if for each object $d\in\cD$, the $\infty$-category $F_{d/}$ defined by the pullback
\[
\begin{tikzcd}
F_{d/}\arrow[d]\arrow[r]\arrow[dr,"\lrcorner",phantom,very near start] & \cD_{d/}\arrow[d] \\
\cC\arrow[r,"F"'] & \cD
\end{tikzcd}
\]
is weakly contractible.
\end{theorem}

\begin{lemma}\label{cofinal_ff}
Let $I,\cC',\cC$ be $\infty$-categories and $j:\cC'\hookrightarrow\cC$ a fully faithful functor. Given a functor $F:I\rightarrow\cC'$, if the composition $jF:I\rightarrow\cC$ is cofinal, then so is $F$.
\end{lemma}

\begin{proof}
Let $c\in\cC'$ and consider the rectangle of $\infty$-categories

\[
\begin{tikzcd}
jF_{j(c)/}\arrow[d]\arrow[r] & \cC'_{c/}\arrow[d]\arrow[r] & \cC_{j(c)/}\arrow[d] \\
I\arrow[r,"F"'] & \cC'\arrow[r,"j"'] & \cC
\end{tikzcd}.
\]
By definition, the whole rectangle is a pullback square. Moreover, because $j$ is fully faithful, the left fibrations $\cC'_{c/}\rightarrow\cC',\cC_{j(c)/}\rightarrow\cC$ have weakly homotopy equivalent fibers so that the right square is a pullback as well. Therefore the left square is also a pullback square. This asserts that $jF_{j(c)/}\simeq F_{c/}$. Because $jF$ is cofinal, Quillen's Theorem A asserts that $jF_{j(c)/}$, and thus $F_{c/}$, are weakly contractible. The same theorem then asserts that $F$ is cofinal.
\end{proof}

\subsection{Localizations}

Given an $\infty$-category $\cC$ and a set of morphisms $W$ in $\cC$, we are interested in formally inverting morphisms in $W$ to obtain a new $\infty$-category $\cC[W^{-1}]$. This is the process of \emph{localizing} $\cC$ at $W$. We use the conventions from \cite[\href{https://kerodon.net/tag/01M4}{Section 6.3}]{kerodon}. 

\begin{definition}
Let $\cC$ be a simplicial set and $W$ a set of $1$-simplices in $\cC$. A morphism $L:\cC\rightarrow\cC[W^{-1}]$ exhibits $\cC[W^{-1}]$ as a \emph{localization} of $\cC$ at $W$ if for any $\infty$-category $\cE$, precomposition by $L$ induces a fully faithful functor $L^*:\Fun(\cC[W^{-1}],\cE)\rightarrow\Fun(\cC,\cE)$ whose essential image consists of morphisms sending elements in $W$ to isomorphisms in $\cE$.
\end{definition}

One benefit of localizations is that they exist for arbitrary pairs $(\cC,W)$ and the assignment $(\cC,W)\mapsto (\cC\rightarrow\cC[W^{-1}])$ can be made functorial. Moreover, localizations are cofinal, see \cite[\href{https://kerodon.net/tag/02N9}{Proposition 7.2.1.10}]{kerodon}. 
Hence, when it comes down to prove that a functor is cofinal, it is enough to prove it is a localization functor. In practice, it is often easier to prove such a statement. As an example, a common strategy is to use \emph{reflective localizations}: for $\cC,\cD$ $\infty$-categories and $\iota:\cD\hookrightarrow\cC$ a fully faithful embedding, if $\iota$ admits a left adjoint, then this left adjoint is a localization functor (\cite[\href{https://kerodon.net/tag/04JL}{Proposition 6.3.3.6}]{kerodon}). 

Another example is from computing colimits. Given a small $\infty$-category $\cC$, a functor $F:\cC^{op}\rightarrow\Cat$ and $\cE\rightarrow\cC$ its unstraightening, there is a functor $\cE\rightarrow\colim{\cC^{op}}{F}$ which is a localization of $\cE$ at all cartesian morphisms (\cite[\href{https://kerodon.net/tag/05S2}{Proposition 7.4.5.1}]{kerodon}). In fact, this localization functor admits an alternative description using the \emph{contravariant refraction} construction (\cite[Definition 7.4.5.8]{kerodon}\footnote{This is made in the covariant case, but as usual it can be adapted to the contravariant case.}). Namely, let $\bar{F}:(\cC^\lhd)^{op}\simeq (\cC^{op})^\rhd\rightarrow\Cat$ be a cocone extending $F$ and $-\infty$ be the initial object of $\cC^{\lhd}$. By \cite[7.4.5.10]{kerodon}, this yields a refraction diagram $\text{Rf}:\cE\rightarrow\bar{F}(-\infty)$ having the property that for each object $u\in\cC$, the composition
\[
\cE_u\hookrightarrow\cE\xrightarrow{\text{Rf}}\bar{F}(-\infty)
\]
is the contravariant transport induced by the unique morphism $-\infty\rightarrow u$. Then, \cite[Theorem 7.4.5.13]{kerodon} states that $\bar{F}:(\cC^{op})^\rhd\rightarrow\Cat$ is a colimit diagram if and only if the functor $\text{Rf}:\cE\rightarrow\bar{F}(-\infty)$ exhibits the latter as the localization of $\cE$ at all cartesian morphisms.

\section{Stratified spaces and exit path categories}\label{section_exit_paths}

Let $T$ be a topological space. We are interested in \emph{stratifications} on $T$, that is a decomposition $T=\bigsqcup\limits_{u\in P}T_u$ into \emph{strata}, with mild compatibility conditions.

One way of formalizing this idea is done in  \cite[Appendix A.5]{ha}. Namely, a stratification on $T$ is simply given by a map $\pi:T\rightarrow P$ where $P$ is a poset and the strata are the preimages $\{T_u:=\pi^{-1}(u)\}_{u\in P}$.

\begin{definition}
Let $P$ be a poset. The \emph{Alexandroff topology} on $P$ is the topology where a subset $U\subset P$ is open if for any $u\in U$ and $v\in P$ such that $u\leq v$, we have $v\in U$.
\end{definition}

In particular, for any $u\in P$, the subsets $P_{u\leq}:=\{v\in P,u\leq v\},P_{u<}:=\{v\in P,u<v\}$ are open. 

\begin{definition}
Let $T$ be a topological space and $P$ a poset. A \emph{$P$-stratification} on $T$ consists of a continuous map $\pi:T\rightarrow P$ where $P$ is endowed with the Alexandroff topology. 
\end{definition}

\begin{definition}
Given a poset $P$, let $P^{\lhd}$ be the poset obtained by adjoining a minimal element $-\infty$ to $P$, that is $P^{\lhd}=\{-\infty\}\sqcup P$ and $-\infty<u$ for any $u\in P$.
\end{definition}

\begin{definition}
Given a topological space $T$, its open cone $C(T)$ is the topological space whose underlying set is $\{\ast\}\sqcup (T\times\R_{>0})$ and $U\subset C(T)$ is open if and only if $U\cap (T\times\R_{>0})$ is open and, if $\ast\in U$, there exists $\epsilon>0$ such that $T\times ]0,\epsilon[\subset U$.
\end{definition}

\begin{definition}
Let $p:T\rightarrow P$ be a stratified space. We define a $P^{\lhd}$-stratification on the cone $C(T)$ by

\[
(x,t)\mapsto
\left\{
\begin{array}{ll}
-\infty & \text{if }t=0. \\
p(x) & \text{otherwise}.
\end{array}
\right.
\]
\end{definition}

For the following definition, first observe that given a topological space $U$ and a stratified space $T\rightarrow P$, the product $U\times T$ is also $P$-stratified by the composition $U\times T\xrightarrow{\pi_T} T\rightarrow P$. Moreover for any $a\in P$, there is an isomorphism of poset $(P_{a<})^{\lhd}\cong P_{a\leq}$.

\begin{definition}
Let $p:T\rightarrow P$ a stratification. We say that this stratification is \emph{conical} if for any point $x\in T$ lying over an element $a\in P$, there exists a $P_{a<}$-stratified space $L$, a topological space $U$ and a commutative square

\[
\begin{tikzcd}
U\times C(L)\arrow[r,hook,"\iota"]\arrow[d] & T\arrow[d] \\
P_{a\leq}\arrow[r,hook] & P 
\end{tikzcd}
\]
where the map $\iota:U\times C(L)\hookrightarrow T$ is an open embedding for which there exists $u\in U$ such that $\iota(u,\ast)=x$
\end{definition}

\begin{definition}
The ordinary category $\Strat$ of conically stratified topological spaces has
\begin{enumerate}
\item Objects: the collection of conically stratified spaces.
\item A morphism $f:(T,P)\rightarrow (S,Q)$ is given by a commutative square of topological spaces
\[
\begin{tikzcd}
T\arrow[r,"f^\sharp"]\arrow[d] & S\arrow[d] \\
P\arrow[r,"f^\flat"'] & Q
\end{tikzcd}.
\]
\end{enumerate}
\end{definition}

Let $\Delta^p=\{(0\leq t_1\leq\cdots\leq t_p\leq 1)\}$ be the geometric $p$-simplex and $[p]=\{0\leq 1\leq\cdots\leq p\}$. One can check that the application $\Delta^p\rightarrow [p]$ sending a point $t=(0\leq t_1\leq\cdots\leq t_p\leq 1\}$ to $\#\{1\leq i\leq p:t_i>0\}$ is continuous and thus defines a $[p]$-stratification on $\Delta^p$. Moreover this map is natural in $[p]$. In particular, for any poset $P$, the induced map $\Top([p],P)\rightarrow\Top(\Delta^p,P)$ extends to a simplicial map $N(P)\rightarrow \Sing(P)$. Here, we use the fact that maps between posets are exactly the continuous maps with respect to their Alexandrov topologies.

This leads to Lurie's model for the \emph{exit path category} of a stratified space.
\begin{proposition}\label{cstrat}\cite[Theorem A.6.4]{ha}
Let $\pi:T\rightarrow P$ be a conically stratified space. Then the pullback of simplicial sets
\[
\begin{tikzcd}
N(P)\times_{\Sing(P)}\Sing(T)\arrow[d]\arrow[r]\arrow[dr, phantom, "\lrcorner", very near start] & \Sing(T)\arrow[d,"\Sing(\pi)"] \\
N(P)\arrow[r] & \Sing(P)
\end{tikzcd}
\]
is an $\infty$-category.
\end{proposition}

\begin{definition}
For a conically stratified space $p:T\rightarrow P$, the $\infty$-category 
\[
N(P)\times_{\Sing(P)}\Sing(T)
\]
is called the exit path category and will be denoted $\exit{P}{T}$. The $\infty$-category $\enter{P}{T}:=\exit{P}{T}^{op}$ is its dual.
\end{definition}

Unwinding the definition, we get that 
\[
\exit{P}{T}_p\cong\Strat((\Delta^p,[p]),(T,P)).
\]
In other word, any $p$-simplex is given by some stratified map
\[
\begin{tikzcd}
\Delta^p\arrow[d]\arrow[r] & T\arrow[d] \\
{[p]}\arrow[r] & P
\end{tikzcd}.
\]

In particular, a $1$-simplex is the datum of a path $\alpha:I\rightarrow T$ where $\alpha(0)\in T_u$ and $\alpha(]0,1])\subset T_v$ for some $u\leq v$ in $P$. This recovers the idea of a path exiting definitely a stratum.

\begin{rem}\label{exit_strict_functor}
Because of \Cref{cstrat}, we get a functor $\exit{-}{-}:\Strat\rightarrow\Cat$. Notice that this functor is strict, in the sense it factors through the ordinary category $\Cat^{\delta}\subset\Cat$ consisting of $\infty$-categories and functors between them. Now let $I$ be an ordinary category. For any functor $F:I\rightarrow\Strat$, we get a strict functor $\exit{-}{-}\circ F:I\rightarrow\Cat$. One can then explicitly model its unstraightening via its \emph{weighted nerve}. This construction will be ubiquitous in this work.
\end{rem}

We now review how to compute mapping spaces in exit path categories.

\begin{definition}\cite[Definition A.7.3]{ha}
Let $\pi:T\rightarrow P$ be a stratified space. Given a chain 
\[
\overrightarrow{u}=(u_0\leq\cdots\leq u_n)\in N(P)_n,
\] 
the simplicial set $\exit{P}{T}[\overrightarrow{u}]$ is defined by the pullback
\[
\begin{tikzcd}
\exit{P}{T}[\overrightarrow{u}]\arrow[d]\arrow[r]\arrow[dr,"\lrcorner",phantom,very near start] & \exit{P}{T}^{\Delta[n]}\arrow[d] \\
\ast\arrow[r,"\overrightarrow{u}"'] & N(P)^{\Delta[n]}
\end{tikzcd}.
\]
\end{definition}

According to \cite[Remark A.7.4]{ha}, this is a Kan complex provided that the stratification is conical. In particular this gives an explicit description for the mapping space between two points $x,y\in T$. 

\begin{lemma}
Let $\pi:T\rightarrow P$ a conical stratification and $x,y\in T$ with $\pi(x)=u, \pi(y)=v$. If $u\not\leq v$, then the mapping space $\exit{P}{T}(x,y)=\emptyset$. Otherwise it is given by the homotopy pullback 
\[
\begin{tikzcd}
\exit{P}{T}(x{, }y)\arrow[d]\arrow[r]\arrow[dr,"\lrcorner",phantom,very near start] & \exit{P}{T}[u\leq v]\arrow[d] \\
\ast\arrow[r,"(x{, }y)"'] & \exit{P}{T}[u]\times\exit{P}{T}[v]
\end{tikzcd}.
\]
\end{lemma}

\begin{proof}
The first case is straight from the definition so assume $u\leq v$. We have by definition that $\exit{P}{T}(x,y)$ is given by the homotopy pullback 
\[
\begin{tikzcd}
\exit{P}{T}(x{, }y)\arrow[d]\arrow[r]\arrow[dr,"\lrcorner",phantom,very near start] & \lf(\exit{P}{T}^{\Delta[1]}\ri)^{\simeq}\arrow[d] \\
\ast\arrow[r,"(x{, }y)"'] & \exit{P}{T}^{\simeq}\times\exit{P}{T}^{\simeq}
\end{tikzcd}.
\]

According to the proof of \cite[Theorem A.7.1]{ha}, we have
\[
\lf(\exit{P}{T}^{\Delta[1]}\ri)^{\simeq}\simeq\bigsqcup\limits_{\overrightarrow{u}\in N(P)_1}\exit{P}{T}[\overrightarrow{u}].
\] 
In particular, restricting to the relevant connected components gives the result.
\end{proof}

To conclude this section, we state a useful result about changes of stratifications on a topological space $T$, see \cite[Theorem 0.3.1]{exodromy}, \cite[Proposition 3.22]{jansen}.

\begin{proposition}
Let $T\rightarrow P$ be a conically stratified space and $P\rightarrow Q$ a morphism between finite posets. If the induced map $T\rightarrow Q$ defines another conical stratification on $T$, then the induced functor 
\[
\exit{P}{T}\rightarrow\exit{Q}{T}
\]
is a localization functor between $\infty$-categories.
\end{proposition}

The relevance of this result is that localizations are cofinal. Hence, in order to prove that a functor $F:\cC\rightarrow\exit{Q}{T}$ is cofinal, it is enough to find a finer stratification $(T,P)$, a functor $\cC\rightarrow\exit{P}{T}$ through which $F$ factors and prove this functor is cofinal. For example, \cite[Theorem A.6.10]{ha} states that whenever $T$ is a locally finite simplicial complex and $P$ is the poset of its simplices, the functor exit $\exit{P}{T}\rightarrow P$ is an equivalence of $\infty$-categories so that the problem amounts to prove cofinality of a functor $\cC\rightarrow P$, which is considerably simpler.  

\subsection{The stratification of $M^n$ by partitions}
We introduce here our first stratification of interest. For the remainder of the section, we fix a closed $d$-dimensional smooth manifold $M$.

\begin{definition}
For any $\us{n}\in\Fin$, the poset $\prt{\us{n}}$ has objects equivalence relations $R\subset\us{n}\times\us{n}$ and partial order $R\leq S$ whenever $S\subset R$. 
\end{definition}

\begin{definition}
Let $\us{n}\in\Fin$. The stratification $\pi^n_M:M^n\rightarrow \prt{\us{n}}$ on $M^n\cong\Map(\us{n},M)$ is defined by sending a point $(x:\us{n}\rightarrow M)$ to the unique equivalence relation $R$ providing a factorization
\[
\begin{tikzcd}
\us{n}\arrow[rr,"x"]\arrow[dr,"f"'] & & M \\
& \us{n}/R\arrow[ur,"x'"',hook] &
\end{tikzcd}.
\]
\end{definition}

\begin{rem}
An equivalence relation on $\us{n}$ is equivalent to a surjective map $\us{n}\twoheadrightarrow\us{k}$ modulo the $\Sigma_k$-action. It will sometimes be convenient to represent such an equivalence relation by a surjective map $\us{n}\rightarrow\us{k}$.
\end{rem}

One observes that the top stratum of the stratification $\pi^n_M:M^n\rightarrow \prt{\us{n}}$ is the configuration space $\Conf(n,M)$. 

Moreover, observe that the assignment $\us{n}\mapsto \prt{\us{n}}$ is contravariantly functorial in $\us{n}$. Namely, given a map $f:\us{n}\rightarrow\us{m}$ , we get a map of posets $f^*:\prt{\us{m}}\rightarrow \prt{\us{n}}$ by sending a surjection $g:\us{m}\twoheadrightarrow\us{k}$ to a surjection $f^*g:\us{n}\twoheadrightarrow\us{k'}$, unique modulo the $\Sigma_k$-action, providing a factorization
\[
\begin{tikzcd}
\us{n}\arrow[d,"f^*g"',two heads]\arrow[r,"f"] & \us{m} \arrow[d,"g",two heads] \\
\us{k'}\arrow[r,hook] & \us{k}.
\end{tikzcd}.
\]

\begin{lemma}\label{functoriality_pi}
For any map $f:\us{n}\rightarrow\us{m}\in\Fin$, the square
\[
\begin{tikzcd}
M^m\arrow[d,"\pi^m_M"']\arrow[r,"f^*"] & M^n\arrow[d,"\pi^n_M"] \\
\prt{\us{m}}\arrow[r,"f^*"'] & \prt{\us{n}} 
\end{tikzcd}
\]
commutes.
\end{lemma}

\begin{proof}
Let $x:\us{m}\rightarrow M$ an element in $M^m$. We get a commutative diagram
\[
\begin{tikzcd}
\us{n}\arrow[d,"f^*\lf(\pi^m_M(x)\ri)"']\arrow[r,"f"] & \us{m}\arrow[d,"\pi_M^m(x)"']\arrow[r,"x"] & M \\
\us{k'}\arrow[r,hook] & \us{k}\arrow[ur,hook] & 
\end{tikzcd}.
\]

Now observe that the compositions $(\us{k'}\hookrightarrow\us{k}\hookrightarrow M)\circ (f^*\lf(\pi_M^m(x)\ri):\us{n}\twoheadrightarrow\us{k'})$ factorizes $f^*x$ as a surjection followed by an injection. By definition, we have $f^*\lf(\pi_M^m(x)\ri)=\pi^n_M(f^*x)$ whence the commutativity of the square.
\end{proof}

A quick check shows it is compatible with the partial order and the composition so that the assignment $\us{n}\mapsto\pi_n:M^n\rightarrow \prt{\us{n}}$ is functorial.

Furthermore, it turns out that the map $\pi^n_M$ is not only continuous but defines a \emph{conical} stratification. 

\begin{lemma}\label{conical_stratification_Mn}
For any $\us{n}\in\Fin$, the stratification $\pi^n_M:M^n\rightarrow \prt{\us{n}}$ is conical.
\end{lemma}

\begin{proof}
Fix a point $x=(x_1,\cdots,x_n)\in M^n$ lying over an element $R\in \prt{\us{n}}$. If $R$ is the maximal element, then $x\in\Conf(n,M)$ so that it admits an open neighbourhood $U\subset\Conf(n,M)$. Then one gets a conical neighbourhood $U\cong U\times C(\emptyset)\hookrightarrow M^n$ over $\prt{\us{n}}_{R\leq}$. Assume now that $R$ is not the maximal element. Using an appropriate chart, we may identify $x$ with a point in $(\R^d)^n$ and it is enough to define a conical neighbourhood in $(\R^d)^n$ around $x$. Let $W$ be the sub vector space of vectors $v=(v_1,\cdots,v_n)\in (\R^d)^n$ satisfying $v_i=v_j$ whenever $x_i=x_j$, and $W^{\perp}$ its orthogonal complement. Define 
\[
\rho:=
\lf\{
\begin{array}{ll}
1 & \text{if }x_1=\cdots=x_n \\
\min_{x_i\neq x_j}\frac{\Norm{x_i-x_j}}{2} & \text{otherwise} 
\end{array}
\ri.
\] 

and $V:=\{v\in W^{\perp},\Norm{v}=\rho\}$. Since $R$ is not maximal, $W^\perp\neq\{0\}$ and thus $V$ is nonempty. Then, for any $v=(v_1,\cdots,v_n)\in V$, we get $x_i\neq x_j$ implies $x_i+v_i\neq x_j+v_j$ and for at least one pair $(i,j)$, $x_i=x_j$ but $x_i+v_i\neq x_j+v_j$. Therefore the $\prt{\us{n}}$-stratification on $x+V$ restricts to a $\prt{\us{n}}_{R<}$-stratification. We then get a $\prt{\us{n}}_{R\leq}$-stratification on the cone $C(x+V):=\{x+tv,t\in [0,1[\}$. Taking an open neighbourhood $U$ around $x$ in the sub vector space $W$, we get the desired conical neighbourhood $U\times C(x+V)\rightarrow  \prt{\us{n}}_{R\leq}$.

\end{proof}

Combining \Cref{functoriality_pi} and \Cref{conical_stratification_Mn}, we deduce a functor 
\[
\pi_M^-:\Fin^{op}\rightarrow\Strat,\us{n}\mapsto (\pi^n_M:M^n\rightarrow \prt{\us{n}}).
\] 

Combining with the enter path functor, we then get a functor
\[
\enter{\prt{\us{-}}}{M^-}:\Fin^{op}\rightarrow\Cat^\delta.
\]

\subsection{The stratification of $\norm{X}$ by nondegenerate simplices}\label{section_stratification_realization}

This is a recollection of \cite[Section 3]{stratified_haya}. For convenience, we will work with    a variant of the simplex category $\Delta^{op}$, the category $\nabla$ of finite linear intervals.

\begin{definition}
The category $\nabla$ has
\begin{enumerate}
\item Objects: finite linearly ordered sets $\la p\ra:=\{\bot<1<\cdots<p<\top\}$ for each integer $p\geq 0$. 
\item Morphisms: a morphism $f:\la p\ra\rightarrow \la q\ra$ is an order preserving function preserving the minimum and maximum elements $\bot,\top$. Compositions and identities are inherited from $\Set$.
\end{enumerate}
\end{definition}

There is an isomorphism of categories $\nabla\cong\Delta^{op}$, see \cite{wraith}. Under this isomorphism, (compositions of) faces $[p]\hookrightarrow [q]$ correspond to surjections $\la q\ra\twoheadrightarrow \la p\ra$ while (compositions of) degeneracies $[p]\twoheadrightarrow [q]$ correspond to injections $\la q\ra\hookrightarrow\la p\ra$.

With this description, we can define the geometric simplex $\Delta^p:=\Map_{\partial}(\la p\ra,I)$ by the topological space of ordered maps $\la p\ra\rightarrow I=[0,1]$ sending the $\bot,\top$ elements to $0,1$ respectively. Thus an element $t\in\Delta^p$ is given by an ordered tuple $(0\leq t_1\leq\cdots\leq t_p\leq 1)$. This definition extends to a functor $\Delta:\nabla^{op}\rightarrow\Top$. Because of this description, using the category $\nabla$ instead of the cosimplicial category $\Delta$ is relevant for our purposes.

\begin{definition}
Let $p\geq 0$. The poset $\pos{p}$ has objects all faces $\theta:\la p\ra\twoheadrightarrow\la q\ra$ out of $\la p\ra$ and partial order given by $\theta\leq\theta'$ if $\theta$ factors through $\theta'$. The stratification $\pi_p:\Delta^p\rightarrow\pos{p}$ is defined by sending an element $t:\la p\ra\rightarrow I$ to the unique face $\theta:\la p\ra\twoheadrightarrow \la q\ra$ providing a factorization $\la p\ra\xrightarrow{\theta}\la q\ra\hookrightarrow I$ of $t$ by a surjection followed by an injection.
\end{definition}

The poset $\pos{p}$ is more classically described as the poset of nondegenerate simplices of the standard simplicial set $\Delta[p]$. With this alternative point of view, we see in particular that we get a functor $\pos{-}:\nabla^{op}\rightarrow\Pos$ and $N(\pos{p})\cong\Sd(\Delta[p])$. Then, the stratification $\pi_p:\Delta^p\rightarrow\pos{p}$ is simply the stratification of $\Delta^p$ by its faces. Moreover this stratification is conical for every $p\geq 0$.

Next we recall the following construction from \cite[Section 3]{stratified_haya}, generalizing the stratification $\pi_p:\Delta^p\rightarrow\pos{p}$. 

\begin{definition}
let $X$ be a simplicial set. The poset $\pos{X}$ has elements the nondegenerate simplices of $X$ and for $\sigma,\tau\in\pos{X}$, a relation $\tau\leq\sigma$ if $\tau$ is a face of $\sigma$. 
\end{definition}

Endowing the posets with the Alexandrov topology, we have a colimit of topological spaces 
\[
\colim{\sigma\in\Delta(X)}{\pos{\sigma}}\cong\pos{X}.
\]
Combining with the geometric realization formula 
\[
\colim{\sigma\in\Delta(X)}{\Delta^\sigma}\cong\norm{X},
\]
this defines a stratification $\pi_X:\norm{X}\rightarrow\pos{X}$. Informally speaking, this stratification sends a point $x\in\norm{X}$ to the maximal nondegenerate simplex containing it.

If the simplicial set $X$ is locally finite, that is each simplex is a face of a finite number of nondegenerate simplices, then the stratification $\pi_X$ is conical (\cite[Theorem 1.1]{stratified_haya}). This gives a diagram $\pi_-:\Delta(X)^\rhd\rightarrow\Strat$ sending a simplex $\sigma\in\Delta(X)$ to $\pi_\sigma:\Delta^\sigma\rightarrow\pos{\sigma}$ and the terminal object $\infty\in\Delta(X)^\rhd$ to $\pi_X:\norm{X}\rightarrow\pos{X}$. Applying the exit path functor, we get a diagram $\text{Exit}\circ\pi_-:\Delta(X)^\rhd\rightarrow\Cat^\delta$. 

\begin{theorem}\label{stratified_realization}(\cite[Theorem 1.2]{stratified_haya})
Let $X$ be a finite simplicial set. Then there is an equivalence of $\infty$-categories
\[
\colim{\sigma\in\Delta(X)}{\exit{\pos{\sigma}}{\Delta^\sigma}}\simeq \exit{\pos{X}}{\norm{X}}.
\]
\end{theorem}

Since taking the opposite category preserves colimits, under the same hypothesis, one deduces the equivalence of $\infty$-categories

\[
\colim{\sigma\in\Delta(X)}{\enter{\pos{\sigma}}{\Delta^\sigma}}\simeq \enter{\pos{X}}{\norm{X}}.
\]

To conclude this section, we explain how the stratification $\pi_X:\norm{X}\rightarrow\pos{X}$ relates to the stratification by partitions. Given a closed manifold $M$, a finite simplicial set $X$, a homeomorphism $h:\norm{X}\xrightarrow{\cong} M$ and $n\geq 0$, we write $\pi^n_{\norm{X}}$ for the stratification $\norm{X}^n\xrightarrow{h^n}M^n\xrightarrow{\pi^n_M}\prt{\us{n}}$, which is also conical.

\begin{proposition}\label{refinement_of_partitions}
For each object $\us{n}\in\Fin$, we have a map between stratified spaces

\[
\begin{tikzcd}
\norm{X^n}\arrow[d,"\pi_{X^n}"']\arrow[r,"\cong"] & \norm{X}^n\arrow[d,"\pi^n_{\norm{X}}"] \\
\pos{X^n}\arrow[r] & \prt{\us{n}}
\end{tikzcd}.
\]
Moreover this map is natural in $\us{n}$.
\end{proposition}

\begin{proof}
We first explicit the map $\pos{X^n}\rightarrow \prt{\us{n}}$. An element in $\pos{X^n}$ is represented by a simplex $\sigma\in X_p^n$ for some $p\geq 0$. This is nothing but a map $\us{n}\rightarrow X_p,i\mapsto\sigma_i\in X_p$. The factorization by a surjection followed by an injection gives a map $\us{n}\twoheadrightarrow\us{n}/R$ which yields an element $R\in \prt{\us{n}}$ onto which $\sigma$ is mapped. Because the degeneracy maps of $X$ are injective, such a factorization does not depend on the simplex representing the element $[\sigma]\in\pos{X^n}$. Moreover, given a face $d:X_p\rightarrow X_q$, we get an induced map $\us{n}/R\rightarrow\us{n}/S$ in the diagram

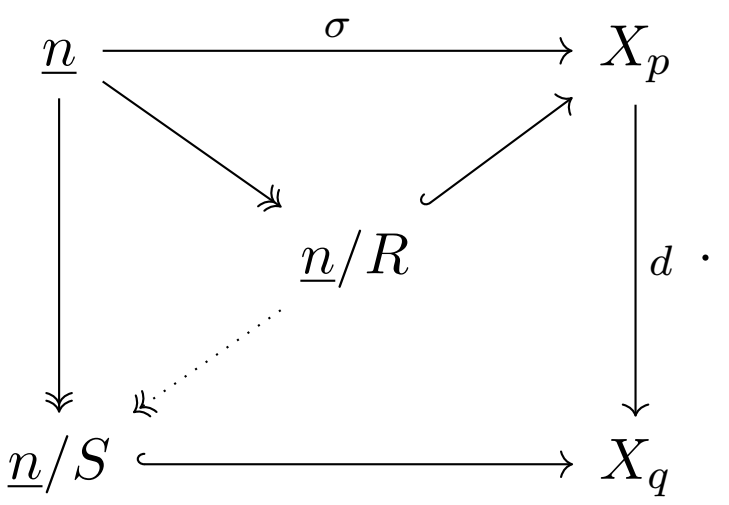
\[
\begin{tikzcd}
\us{n}\arrow[rr,"\sigma"]\arrow[dr,two heads]\arrow[dd,two heads] & & X_p\arrow[dd,"d"] \\
& \us{n}/R\arrow[ur,hook]\arrow[dl,dotted,two heads] & \\
\us{n}/S\arrow[rr,hook] & & X_q
\end{tikzcd}.
\]
This ensures that the map is compatible with the partial order. The map $\pos{X^n}\rightarrow \prt{\us{n}}$ is thus well defined.
Notice also that this morphism is natural in $\us{n}$.

Now let $x\in\norm{X^n}$. This element is represented by an element $t\in\Delta^\sigma$ for some $\sigma\in X^n_p$. Observe that $t$ factors as $\la p\ra\xrightarrow{f}\la r\ra\xrightarrow{t'} I$ where $f$ is a surjection and $t'$ an injection. Thus $t$ represents the same element in $\norm{X^n}$ as $(t'\in\Delta^{f\sigma})$. The upshot is that we may assume that $t:\la p\ra\rightarrow I$ is injective.

Furthermore, the upper map sends $x$ to the tuple $(x_i\in\norm{X})_{i\in\us{n}}$ where $x_i$ is represented by $(t\in\Delta^{\sigma_i})$. The right vertical map sends this tuple to the equivalence relation $S$ defined by $i\sim_S j$ if and only if $x_i=x_j$. It remains to show that $R=S$. By definition of the map $\pos{X^n}\rightarrow \prt{\us{n}}$, the lower right composition maps $x$ to the equivalence relation $R$ given by $i\sim_R j$ if and only if $\sigma_i=\sigma_j$. If $i\sim_R j$, then $\sigma_i=\sigma_j$ thus $(t\in\Delta^{\sigma_i})$ and $(t\in\Delta^{\sigma_j})$ represent the same element. Therefore $x_i=x_j$ whence $i\sim_S j$. Conversely, assume that $i\sim_S j$ so that $x_i=x_j$. The Eilenberg-Zilberg lemma provides unique nondegenerate simplices $ez(\sigma_i),ez(\sigma_j)$ for which $\sigma_i,\sigma_j$ are degeneracies respectively given by $s_i:\la p_i\ra\hookrightarrow\la p\ra,s_j:\la p_j\ra\hookrightarrow\la p\ra$. These degeneracies identify $(t\in\Delta^{\sigma_i})$ and $(t\in\Delta^{\sigma_j})$ with $(ts_i\in\Delta^{p_i}_{ez(\sigma_i)}),(ts_j\in\Delta^{p_j}_{ez(\sigma_j)})$. The element $ts_i$ gives an injective map $\la p_i\ra\xrightarrow{s_i}\la p\ra\xrightarrow{t} I$ hence is not in the image of any face inclusion $\Delta^{p_i-1}_{d(ez(\sigma_i))}\hookrightarrow\Delta^{ez(\sigma_i)}$ and likewise for $ts_j$. Because $ts_i$ and $ts_j$ represent the same element, we necessarily have $ez(\sigma_i)=ez(\sigma_j)$ and $(ts_i\in\Delta^{ez(\sigma_i)})=(ts_j\in\Delta^{ez(\sigma_j)})$. Furthermore, because $t:\la p\ra\rightarrow I$ is injective, we must have $s_i=s_j:\la p_i\ra\rightarrow\la p\ra$. This implies that $\sigma_i=\sigma_j$ whence $i\sim_R j$. 

Finally, because of the naturality in $\us{n}$ of the map $\norm{X^n}\rightarrow\norm{X}^n$ and the map $\pos{X^n}\rightarrow \prt{\us{n}}$, we get the naturality of the stratified map.

\end{proof}

\section{The configuration category of a manifold}\label{section_conf_cat}

Let $M$ be a smooth $d$-dimensional closed manifold.  We translate the particle model defined in \cite[Section 3.1]{confcat} in the setting of quasicategories, using the straightening-unstraightening theorem.  

Consider the functor 
\[
\enter{\prt{\us{-}}}{M^-}:\Fin^{op}\rightarrow\Cat^{\delta}.
\]
By \Cref{exit_strict_functor}, its unstraightening can be modeled as a weighted nerve.

\begin{definition}
The cartesian fibration $\X_M\rightarrow\Fin$ is the weighted nerve of the functor 
\[
\text{Enter}^{\prt{\us{-}}}(M^-):\Fin^{op}\rightarrow\Cat^{\delta}.
\]
\end{definition}

We observe that the objects of $\X_M$ are given by pairs $(\us{n}\in\Fin,x\in M^n)$. 

\begin{definition}
The \emph{configuration category} of $M$, $\con(M)$, is the full subcategory of $\X_M$ where objects are pairs $(\us{n},x\in M^n)$ for which $x\in \Conf(n,M)$.
\end{definition}

We devote the remaining of this section to prove two relevant properties of $\con(M)$. We first prove that there is a localization functor $\X_M\rightarrow\con(M)$ (\Cref{con(M)_localization}).

In addition, the $\infty$-categories $\con(M),\cDisk_{/M}$ come with functors $\con(M)\rightarrow\Fin,\cDisk_{/M}\rightarrow\Fin$. One then proves (\Cref{equiv_disks}) an equivalence in $(\Cat)_{/\Fin}$
\[
\con(M)\simeq\cDisk_{/M}.
\]

\subsection{The configuration category as a localization}

Let $\cart(\Fin^{surj})$ be the set of cartesian lifts of surjective maps $f:\us{n}\twoheadrightarrow\us{k}$ relative to the cartesian fibration $\X_M\rightarrow\Fin$.

\begin{proposition}\label{con(M)_localization}
The inclusion $\iota:\con(M)\hookrightarrow\X_M$ admits a left adjoint $L:\X_M\rightarrow \con(M)$.
\end{proposition}

\begin{lemma}\cite[Lemma 2.10]{volpe}
Let $(S,P)\rightarrow (T,Q)$ be a map between conically stratified spaces such that
\begin{enumerate}
\item the induced map $P\rightarrow Q$ is injective and full,
\item it induces a homeomorphism $S\cong P\times_Q T$.
\end{enumerate}
Then the functor $\exit{P}{S}\rightarrow\exit{Q}{T}$ is fully faithful.
\end{lemma}

\begin{lemma}\label{f_surj_fully_faithful}
Let $f:\us{n}\twoheadrightarrow\us{m}$ be a surjective map. Then the induced functor $\enter{\prt{\us{m}}}{M^m}\rightarrow\enter{\prt{\us{n}}}{M^n}$ is fully faithful.
\end{lemma}

\begin{proof}
It follows from the previous lemma upon verifying that the map $\prt{\us{m}}\rightarrow \prt{\us{n}}$  is injective and full, and that the square
\[
\begin{tikzcd}
M^m\arrow[d]\arrow[r,"f^*"] & M^n\arrow[d] \\
\prt{\us{m}}\arrow[r,"f^*"'] & \prt{\us{n}}
\end{tikzcd}
\]
is a pullback square of topological spaces.

For the first claim, let $R,R'\in \prt{\us{m}}$ such that $f^*(R)=f^*(R')\in \prt{\us{n}}$. Let $(a,b)\in R$. Because $f$ is surjective, there exists $\bar{a},\bar{b}\in\us{n}$ such that $f(\bar{a})=a,f(\bar{b})=b$. Thus, by definition, $(\bar{a},\bar{b})\in f^*(R)$. Therefore $(\bar{a},\bar{b})\in f^*(R')$ and, again by definition, we necessarily have $(a,b)\in R'$. This implies that $R\subset R'$ and by symmetry, $R=R'$. A similar argument shows that if $f^*(R)\subset f^*(R')$, then $R\subset R'$ whence the fullness. 

For the second claim, observe that the map $M^m\rightarrow M^n$ is an embedding whose image consists of all points $x:\us{n}\rightarrow M$ factoring through the surjection $f:\us{n}\twoheadrightarrow\us{m}$. By inspection, the image is the fiber product $\prt{\us{m}}\times_{\prt{\us{n}}}M^n$ and thus the square is indeed a pullback square.
\end{proof}

\begin{lemma}
The objects of $\con(M)$ are $\cart(\Fin^{surj})$-local.
\end{lemma}

\begin{proof}
Let $c\in \con(M)$ and $\alpha:x\rightarrow y$ a cartesian lift of a surjective morphism $f:\us{n'}\twoheadrightarrow\us{n}$. Without loss of generality, we identify $x$ with $f^*y$ and the induced path $f^*y\rightarrow f^*y$ in $\exit{\prt{\us{n'}}}{M^{n'}}$ is the constant one. Then the induced map $\X_M(y,c)\rightarrow\X_M(f^*y,c)$ is weakly homotopy equivalent to the map
\[
\bigsqcup\limits_{g:\us{n}\rightarrow\us{m}}\text{Enter}^{\prt{\us{n}}}(M^n)(y,g^*c)\xrightarrow{f^*}\bigsqcup\limits_{g':\us{n'}\rightarrow\us{m}}\text{Enter}^{\prt{\us{n'}}}(M^{n'})(f^*y,(g')^*c)
\]
and it is enough to prove the weak homotopy equivalence componentwise. Fix a map $g':\us{n'}\rightarrow\us{m}$. We distinguish two cases:
\begin{enumerate}
\item $g'$ factors through $f$: write $g'=fg$ with $g:\us{n}\rightarrow\us{m}$. Note that $g$ is unique because $f$ is surjective. \Cref{f_surj_fully_faithful} asserts that the functor $f^*:\enter{\prt{\us{n}}}{M^n}\rightarrow\enter{\prt{\us{n'}}}{M^{n'}}$ is fully faithful so that we get a weak homotopy equivalence 

\[
\text{Enter}^{\prt{\us{n}}}(M^n)(y,g^*c)\rightarrow\text{Enter}^{P(\us{n')}}(M^{n'})(f^*y,f^*g^*c).
\]

\item $g'$ does not factor through $f$: let $R\in \prt{\us{n'}}$ the stratum where $f^*y$ lies and $R(f)\in \prt{\us{n'}}$ the relation induced by $f$, that is $(a,b)\in R(f)$ if and only if $f(a)=f(b)$. By definition, the map $f^*y:\us{n}\rightarrow M$ factors through $f$ which implies $R\leq R(f)$. Now let $S$ be the stratum where $(g')^*c$ lies. Because $g'$ does not factor through $f$ and $c$ is injective, we have $S\not\leq R(f)$. In particular $S\not\leq R$ and therefore $\enter{\prt{\us{n'}}}{M^{n'}}(f^*y,(g')^*c)=\emptyset$.
\end{enumerate}

To summarize, we have proved that each component on the right side is either empty, or is weakly homotopy equivalent to exactly one term on the left side. Therefore the map is indeed a weak homotopy equivalence.
\end{proof}

\begin{proof}[Proof of \Cref{con(M)_localization}]
By \cite[Proposition 5.1.10]{land}, it is enough to specify for each object $(x:\us{n}\rightarrow M)\in\X_M$, an object $L(x)\in \con(M)$ and a morphism $\alpha_x: x\rightarrow\iota(L(x))$ for which the composition

\begin{equation}\tag{\ensuremath{\ast}}\label{unit_localization}
\con(M)(L(x),c)\rightarrow\X_M(\iota(L(x)),\iota(c))\xrightarrow{\alpha_x^*}\X_M(x,\iota(c))
\end{equation}

is a weak homotopy equivalence. Consider the factorization
\[
\begin{tikzcd}
\us{n}\arrow[dr,"x"']\arrow[rr,"f",two heads] & & \us{n'}\arrow[dl,"L(x)",hook] \\
& M & 
\end{tikzcd}
\]
by a surjection followed by an injection. Define $\alpha_x:x\rightarrow L(x)$ to be the cartesian lift of $f$ with codomain $L(x)$ relatively to the cartesian fibration $\X_M\rightarrow\Fin$. Then $\alpha_x\in \cart(\Fin^{surj})$ and $L(x)\in \con(M)$.

The map $\con(M)(L(x),c)\rightarrow\X_M(\iota(L(x)),\iota(c))$ is a weak homotopy equivalence because the inclusion $\iota:\con(M)\hookrightarrow\X_M$ is fully faithful. Moreover, the map $\X_M(\iota(L(x)),\iota(c))\xrightarrow{\alpha_x^*}\X_M(x,\iota(c))$ is also a weak homotopy equivalence because $L(x)$ is $\cart(\Fin^{surj})$-local. Therefore the map \eqref{unit_localization} is also a weak homotopy equivalence.
\end{proof}

\subsection{The equivalence with $\cDisk_{/M}$}

Let $\Man_d$ be the ordinary category of $d$-dimensional smooth manifolds without boundary and embeddings between them.
Consider the topological category with the same objects and for two manifolds $M,N$, the enrichment given by $\Emb(M,N)$, the space of embeddings $M\hookrightarrow N$ endowed with the weak $C^\infty$-topology. Its coherent nerve defines an $\infty$-category written $\cMan_d$. From now on, the dimension will stay fixed hence we drop it from the notation. We have an obvious functor $\Man\rightarrow\cMan$.

Let $\Disk\subset\Man,\cDisk\subset\cMan$ be the full subcategories whose objects are manifolds of the form $\bigsqcup\limits_{n}\R^d$ for some $n\geq 0$. For a manifold $M\in\Man$, one then gets the overcategories $\Disk_{/M},\cDisk_{/M}$. 

Consider the path component functor $\pi_0:\cDisk\rightarrow\Fin$. By composition, we get a functor $\cDisk_{/M}\rightarrow\cDisk\rightarrow\Fin$ we also write $\pi_0$.

\begin{proposition}\label{equiv_disks}
There is an equivalence in $(\Cat)_{/\Fin}$
\[
\con(M)\simeq\cDisk_{/M}.
\]

\end{proposition}

This equivalence has already been proved in \cite[Section 3.2]{confcat} in the setting of (fiberwise) complete Segal spaces. However, the proof of \Cref{equiv_disks} does not translate easily to our quasicategorical setting. Instead, our proof relies heavily on the proof of \cite[Proposition 2.19]{factorization}, exhibiting the functor $\Disk_{/M}\rightarrow\cDisk_{/M}$ as a localization at isotopies.

\begin{rem}
While the results in \cite[Section 2]{factorization} are done for \emph{topological} manifolds, it is pointed out in \cite[Remark 3.29]{factorization} that the arguments adapt verbatim for \emph{smooth} manifolds, see also \cite[Appendix A]{arakawa}.
\end{rem}

\begin{lemma}\label{core_con(M)}
There is a weak homotopy equivalence 
\[
\con(M)^{\simeq}\simeq\bigsqcup\limits_{n\geq 0}\Conf(n,M)_{h\Sigma_n}
\]
where each $\Sigma_n$ acts by permutation of points.
\end{lemma}

\begin{proof}
By \Cref{core_cartesian_fib}, the map $(\X_M)^{\simeq}\rightarrow\Fin^\simeq$ is the unstraightening of the functor 
\[
\text{Enter}^{\prt{\us{-}}}(M^-)^{\simeq}:\Fin^\simeq\rightarrow\cS.
\] 
We have a weak homotopy equivalence 
\[
\Fin^{\simeq}\simeq\bigsqcup\limits_{n\geq 0}B\Sigma_n
\] 
and for each $n\geq 0$, the $\Sigma_n$-action on $\enter{\prt{\us{n}}}{M^n}^{\simeq}$ is given by permutation of points. For each $n\geq 0$, define $((\X_M)^\simeq)_n$ to be the pullback
\[
\begin{tikzcd}
((\X_M)^\simeq)_n\arrow[d]\arrow[r]\arrow[dr,"\lrcorner",phantom,very near start] & (\X_M)^\simeq\arrow[d] \\
B\Sigma_n\arrow[r,hook] & \bigsqcup\limits_{k\geq 0}B\Sigma_k
\end{tikzcd}.
\]

One then gets a homotopy fiber sequence
\[
\text{Enter}^{\prt{\us{n}}}(M^n)^{\simeq}\rightarrow((\X_M)^{\simeq})_n\rightarrow B\Sigma_n
\]
exhibiting $((\X_M)^{\simeq})_n$ as the space of coinvariants $\enter{\prt{\us{n}}}{M^n}^{\simeq}_{h\Sigma_n}$.

Moreover we have 
\[
\text{Enter}^{\prt{\us{n}}}(M^n)^\simeq \simeq \bigsqcup\limits_{R\in \prt{\us{n}}}M^n[R]
\] where $M^n[R]$ is the $R$-stratum of the stratification. Because of this decomposition, the previous homotopy fiber sequence restricts to another homotopy fiber sequence
\[
\Conf(n,M)\rightarrow (\con(M)^{\simeq})_n\rightarrow B\Sigma_n
\]
exhibiting $(\con(M)^{\simeq})_n$ as the space of coinvariants $\Conf(n,M)_{h\Sigma_n}$. Assembling all the path components yields the weak homotopy equivalence
\[
\con(M)^{\simeq}\simeq\bigsqcup\limits_{n\geq 0}\Conf(n,M)_{h\Sigma_n}.
\]
\end{proof}

Since embeddings between manifolds are injective and by definition of the configuration category using the weighted nerve, one gets that the assignment $N\mapsto\con(N)$ extends to a strict functor $\con(-):\Man\rightarrow\Cat^\delta$. In particular it gives another functor $\con(-):\Disk_{/M}\rightarrow\Cat^\delta$. 

\begin{definition}
The $\infty$-category $P_M$ is defined to be the weighted nerve of the functor $\con(-):\Disk_{/M}\rightarrow\Cat^\delta$. The full subcategory $Q_M\subset P_M$ has for objects pairs $(e:D\hookrightarrow M,c:\us{n}\hookrightarrow D)$ for which $c:\us{n}\hookrightarrow D$ induces a bijection on path components. 
\end{definition}

\begin{lemma}\label{equiv_Q_M_Disk_M}
The induced functor $Q_M\rightarrow\Disk_{/M}$ is an equivalence of $\infty$-categories.
\end{lemma}

\begin{proof}
For any embedding $(D\hookrightarrow M)\in\Disk_{/M}$, taking the configuration $c:\us{n}\hookrightarrow D$ consisting of the centers of the disks gives an element $(D,c)\in Q_M$ being sent to $D$. Thus the functor is essentially surjective. 

Now let $(D,c),(D',c')\in Q_M$. Because the inclusion $Q_M\hookrightarrow P_M$ is fully faithful and the functor $P_M\rightarrow\Disk_{/M}$ is a cocartesian fibration, for any embedding $\eta:(D\hookrightarrow M)\hookrightarrow (D'\hookrightarrow M)$ we get a homotopy fiber sequence
\[
\con(D')(\eta c,c')\rightarrow Q_M((D,c),(D',c'))\rightarrow\Disk_{/M}.
\]  

Looking at each connected component of $D'$ separately, we get that 

\[
\con(D')(\eta c,c')\simeq\prod_{i\in\us{m}}\con(\{i\}\times\R^d)(\pi_0(\eta c)^{-1}(i),c'_i) 
\]

where each $c'_i$ consists of the single point in the component $\{i\}\times\R^d\subset D'$. By \cite[Lemma 3.4]{confcat}, $\con(\{i\}\times\R^d)(\pi_0(\eta c)^{-1}(i),c'_i)$ is weakly contractible hence $\con(D')(\eta c,c')$ is also weakly contractible. This implies that the map $Q_M((D,c),(D',c'))\rightarrow\Disk_{/M}(D,D')$ is a weak homotopy equivalence showing that the functor is fully faithful.

\end{proof}

\begin{definition}
Given a point $x\in M^n$, an element $R\in \prt{\us{n}}$ and a stratum $M^n[R]\subset M^n$, the \emph{homotopy link} $\text{Holink}(x,M^n,M^n[R])$ is defined by the pullback
\[
\begin{tikzcd}
\text{Holink}(x,M^n,M^n[R])\arrow[d]\arrow[r]\arrow[dr,very near start,phantom,"\lrcorner"] & \exit{\prt{\us{n}}}{M^n}_{x/}\arrow[d,"\text{ev}_1"] \\
\Sing(M^n[R])\arrow[r,hook] & \exit{\prt{\us{n}}}{M^n} 
\end{tikzcd}.
\]
\end{definition}

Unwinding the definition, it is the space of paths $\alpha:I\rightarrow M^n$ such that $\alpha(0)=x$ and $\alpha(]0,1])\subset M^n[R]$. Moreover, taking maximal subgroupoids and using the fact that 
\[
\exit{\prt{\us{n}}}{M^n}^{\simeq}\simeq\bigsqcup\limits_{R\in \prt{\us{n}}}\Sing(M^n[R]),
\] 
we deduce that $\lf(\exit{\prt{\us{n}}}{M^n}_{x/}\ri)^{\simeq}\simeq\bigsqcup\limits_{R\in \prt{\us{n}}}\text{Holink}(x,M^n,M^n[R])$.

Following \cite[Proposition 3.3]{confcat}, Boavida de Brito and Weiss observed that the homotopy link is a local concept: given an open neighbourhood $U$ of $x$, restricting the paths to $U$ provides a weak homotopy equivalence $\text{Holink}(x,U,M^n[R]\cap U)\hookrightarrow \text{Holink}(x,M,M^n[R])$. Notice that this inclusion is compatible with the $\Sigma_n$-action by permutation of points, so that it induces another weak homotopy equivalence 
\begin{align*}
\lf(\bigsqcup\limits_{R\in \prt{\us{n}}}\text{Holink}(x,U,M^n[R]\cap U)\ri)_{h\Sigma_n}
\hookrightarrow \lf(\bigsqcup\limits_{R\in \prt{\us{n}}}\text{Holink}(x,M,M^n[R])\ri)_{h\Sigma_n}.
\end{align*}

Moreover, for any $f:\us{m}\rightarrow\us{n}$, write $\Sigma_f:=\Sigma_{f^{-1}(1)}\times\cdots\times\Sigma_{f^{-1}(n)}$. Notice that $\Sigma_f$ is a subgroup of $\Sigma_m$ so that the previous equivariant inclusion restricts to a $\Sigma_f$ equivariant inclusion.

\begin{lemma}\label{local_holink}
Let $c\in\con(M)$ and $U$ a disk embedded in $M$ containing $c$. Then the embedding $U\hookrightarrow M$ induces a weak homotopy equivalence
\[
\lf(\con(U)_{/c}\ri)^{\simeq}\xrightarrow{\simeq}\lf(\con(M)_{/c}\ri)^{\simeq}.
\]
\end{lemma}

\begin{proof}
Write $c:\us{n}\hookrightarrow M$.
By the previous observation on homotopy links, the result will follow upon proving the weak homotopy equivalence

\[
(\con(M)_{/c})^{\simeq}\simeq\bigsqcup\limits_{f\in\Fin(\us{m},\us{n})_{\Sigma_m}}\text{Holink}(f^*c,\Conf(m,M))_{h\Sigma_f}
\]

and likewise for $U$ (remember that $\Conf(m,M)$ is the top stratum of the stratified space $(M^m,\prt{\us{m}})$. 






The cartesian fibration $\X_M\rightarrow\Fin$ induces for any $(f:\us{m}\rightarrow\us{n})\in\Fin_{/\us{n}}$ a fiber sequence of $\infty$-categories
\[
\text{Enter}^{\prt{\us{m}}}(M^m)_{/f^*c}\rightarrow(\X_M)_{/c}\rightarrow\Fin_{/\us{n}}.
\]

Taking cores gives a homotopy fiber sequence 

\[
\lf(\text{Enter}^{\prt{\us{m}}}(M^m)_{/f^*c}\ri)^{\simeq}\rightarrow\lf((\X_M)_{/c}\ri)^{\simeq}\rightarrow\lf(\Fin_{/\us{n}}\ri)^{\simeq}.
\]

Now observe that the assignment $(f:\us{m}\rightarrow\us{n})\mapsto (\us{\norm{f^{-1}(1)}},\cdots,\us{\norm{f^{-1}(n)}})$ defines an equivalence of (ordinary) categories $\Fin_{/\us{n}}\simeq\Fin^n$. Indeed, it is essentially surjective and the bijection on hom sets comes from the fact that specifying a map $g:\us{m}\rightarrow\us{m'}$ above two maps $f:\us{m}\rightarrow\us{n},f':\us{m'}\rightarrow\us{n}$ amounts to specify maps $g_i:f^{-1}(i)\rightarrow(f')^{-1}(i)$ on all fibers for $i=1,\cdots,n$.

Because of this equivalence, we have 
\begin{align*}
(\Fin_{/\us{n}})^{\simeq} & \simeq (\Fin^n)^{\simeq} \\
& \simeq\bigsqcup\limits_{m_1,\cdots,m_n,\geq 0}B\Sigma_{m_1}\times\cdots B\Sigma_{m_n} \\
& \simeq \bigsqcup\limits_{f\in \Fin(\us{m},\us{n})_{\Sigma_m}}B\Sigma_f.
\end{align*}
Therefore, for a given $f:\us{m}\rightarrow\us{n}$, by restricting to an appropriate subspace $\lf((\X_M)_{/c}\ri)^{\simeq}_f\subset \lf((\X_M)_{/c}\ri)^{\simeq}$,  the homotopy fiber sequence over $f$ becomes
\[
\bigsqcup\limits_{R\in \prt{\us{m}}}\text{Holink}(f^*c,M^m[R])\rightarrow\lf((\X_M)_{/c}\ri)^{\simeq}_f\rightarrow B\Sigma_f.
\]

The group $\Sigma_f$ acts of fibers by permutation of components, hence this homotopy fiber sequence is identified with the Borel construction

\begin{align*}
\bigsqcup\limits_{R\in \prt{\us{m}}}\text{Holink}(f^*c,M^m[R]) & \rightarrow \lf(\bigsqcup\limits_{R\in \prt{\us{m}}}\text{Holink}(f^*c,M^m[R])\ri)_{h\Sigma_f} \\
& \rightarrow B\Sigma_f.
\end{align*}

Assembling all path components, we deduce that 
\[
((\X_M)_{/c})^{\simeq}\simeq \bigsqcup\limits_{f\in\Fin(\us{m},\us{n})_{\Sigma_m}}\lf(\bigsqcup\limits_{R\in \prt{\us{m}}}\text{Holink}(f^*c,M^m[R])\ri)_{h\Sigma_f}.
\]

Note that for each $f:\us{m}\rightarrow\us{n}$, the $\Sigma_f$-action on $M^m$ preserves the top stratum $\Conf(m,M)$. Indeed, permuting the labels of a configuration $(x_1,\cdots,x_m)$ produces another configuration. Therefore the inclusion

\[
\text{Holink}(f^*c,\Conf(m,M))\hookrightarrow \bigsqcup\limits_{R\in \prt{\us{m}}}\text{Holink}(f^*c,M^m[R])
\]
provides a well-defined map

\[
\text{Holink}(f^*c,\Conf(m,M))_{h\Sigma_f}\rightarrow \lf(\bigsqcup\limits_{R\in \prt{\us{m}}}\text{Holink}(f^*c,M^m[R])\ri)_{h\Sigma_f}.
\]

By definition of $\con(M)$, one finally deduces the weak homotopy equivalence
\[
((\con(M)_{/c})^{\simeq}\simeq \bigsqcup\limits_{f\in\Fin(\us{m},\us{n})_{\Sigma_m}}\text{Holink}(f^*c,\Conf(m,M))_{h\Sigma_f}.
\]

\end{proof}

\begin{proof}[Proof of \Cref{equiv_disks}]
The functor $\con(-):\Disk_{/M}\rightarrow\Cat$ extends to a cocone $\Disk_{/M}^{\rhd}\rightarrow\Cat$ whose terminal object is sent to $\con(M)$. This defines a refraction functor $P_M\rightarrow\con(M)$ sending an object $(e:D\hookrightarrow M,c\in\con(D))$ to the configuration $e(c)\in\con(M)$. The restriction to $Q_M$ fits into a diagram

\[
\begin{tikzcd}
\Disk_{/M}\arrow[dr,"\pi_0"'] & Q_M\arrow[d]\arrow[l]\arrow[r] & \con(M)\arrow[dl] \\
& \Fin &
\end{tikzcd}
\]

which is commutative by definition of $Q_M$.

We claim that the restriction $Q_M\rightarrow\con(M)$ exhibits $\con(M)$ as a localization of $Q_M$ at $J_M$, where $J_M$ is the set of morphisms sent to isotopies under the equivalence $Q_M\xrightarrow{\simeq}\Disk_{/M}$ from \Cref{equiv_Q_M_Disk_M}. Given a morphism $(D,c)\rightarrow (D',c')$ in $J_M$, because the embedding $D\hookrightarrow D'$ is an isotopy and both configurations $c:\us{n}\hookrightarrow D,c':\us{n'}\hookrightarrow D'$ are bijective on path components, the induced morphism $(c\rightarrow c')\in\con(D')$ is necessarily an isomorphism (both configurations have exactly one point on each path component of $D'$). Thus its image in $\con(M)$ is an isomorphism as well. This provides a functor $Q_M[J_M^{-1}]\rightarrow\con(M)$. The equivalence $Q_M\xrightarrow{\simeq}\Disk_{/M}$ yields an equivalence $Q_M[J_M^{-1}]\xrightarrow{\simeq}\Disk_{/M}[\text{isotopies}^{-1}]$ and the latter is equivalent to $\cDisk_{/M}$ by \cite[Proposition 2.19]{factorization} and \cite[Remark A.8]{arakawa}. Hence it is enough to prove that the functor $Q_M[J_M^{-1}]\rightarrow\con(M)$ is an equivalence of $\infty$-categories as well. The proof of this statement follows closely the proof of \cite[Proposition 2.19]{factorization}, checking the equivalence on cores and spaces of morphisms.

By \cite[Lemma 2.12]{factorization}, there is a weak homotopy equivalence
\[
(\cDisk_{/M})^{\simeq}\simeq\bigsqcup\limits_{n\geq 0}\Conf(n,M)_{h\Sigma_n}
\]
over $\Fin^\simeq$, given by quotienting the action by permutation of components. Combining with the weak homotopy equivalence $(Q_M[J_M^{-1}])^{\simeq}\simeq (\cDisk_{/M})^{\simeq}$ and \Cref{core_con(M)}, we deduce that the map 
\[
(Q_M[J_M^{-1}])^{\simeq}\rightarrow(\con(M))^{\simeq}
\] 
is a weak homotopy equivalence over $\Fin^\simeq$.

Now given an object $(D\hookrightarrow M)\in\Disk_{/M}$, consider the ladder

\[
\begin{tikzcd}
Q_D[J_D^{-1}]^{\simeq}\arrow[r]\arrow[d] & \con(D)^{\simeq} \arrow[d] \\
\Fun([1],Q_M[J_M^{-1}])^{\simeq}\arrow[r]\arrow[d,"\text{ev}_1"'] & \Fun([1],\con(M))^{\simeq}\arrow[d,"\text{ev}_1"] \\
Q_M[J_M^{-1}]\arrow[r] & \con(M)^{\simeq}
\end{tikzcd}.
\]

Again by \cite[Proposition 2.19]{factorization}, the left vertical sequence is a homotopy fiber sequence. Moreover consider a configuration $(c:\us{n}\hookrightarrow M)\in\con(M)$ restricting to a homotopy equivalence $(c:\us{n}\xrightarrow{\simeq} D)\in\con(D)$. The fiber of the map $\text{ev}_1:\Fun([1],\con(M))^{\simeq}\rightarrow\con(M)^{\simeq}$ over $c$ is given by $(\con(M)_{/c})^{\simeq}$.  By \Cref{local_holink}, this is equivalent to $(\con(D)_{/c})^{\simeq}$. Because the configuration $c$ consists of exactly one point per path component of $D$, \cite[Lemma 3.4]{confcat} tells us that $c$ is terminal in $\con(D)$ hence there is a further weak homotopy equivalence $\con(D)^{\simeq}\simeq (\con(D)_{/c})^{\simeq}$. The upshot is that the right vertical sequence of the ladder is also a homotopy fiber sequence.

Because we have a map between two homotopy fiber sequences whose top and bottom maps are weak homotopy equivalences, so is the middle one. 

To summarize, we have a commutative diagram over $\Fin$ 

\[
\begin{tikzcd}
\Disk_{/M}\arrow[d] & Q_M\arrow[l,"\simeq"']\arrow[d]\arrow[dr] & \\
\cDisk_{/M} & Q_M[J_M^{-1}]\arrow[l,"\simeq"]\arrow[r,"\simeq"'] & \con(M)
\end{tikzcd}
\]

showing the equivalence $\cDisk_{/M}\simeq\con(M)$ in $(\Cat)_{/\Fin}$.

\end{proof}

\section{A variant of the configuration category for the interval}\label{section_conf_cat_interval}

Recall the functor 
\[
(\pi_\bullet:\Delta^\bullet\rightarrow\pos{\bullet}):\nabla^{op}\rightarrow\Strat
\] 
from \Cref{section_stratification_realization}. By composing with the enter path functor $\enter{-}{-}:\Strat\rightarrow\Cat^{\delta}$, we get a functor 
\[
\enter{\pos{\bullet}}{\Delta^\bullet}:\nabla^{op}\rightarrow\Cat^{\delta}.
\] 
By \Cref{exit_strict_functor}, the unstraightening of this functor can be modeled as a weighted nerve.

\begin{definition}
The functor $\F\rightarrow\nabla$ is the weighted nerve of the functor $\enter{\pos{\bullet}}{\Delta^\bullet}:\nabla^{op}\rightarrow\Cat^{\delta}$.
\end{definition}

By inspection, objects of $\F$ are given by a pair $(\la p\ra\in\nabla,u\in\Delta^p)$. This provides the following definition similar to the definition of configuration category.

\begin{definition}
We define the full subcategory $\con^<(I)\subset\F$ to have objects  pairs $(\la p\ra\in\nabla,u\in\Delta^p)$ where $u:\la p\ra\hookrightarrow I$ is an injective map.
\end{definition}

As for the configuration category of a manifold, we prove that there is a localization functor $\F\rightarrow\con^<(I)$ and an equivalence of $\infty$-categories $\con^<(I)\simeq\nabla$. The proof of the first claim follows exactly the proof of \Cref{con(M)_localization}.

\begin{definition}
We define $\cart(\nabla^{surj})$ to be the set of cartesian lifts of face maps $f:\la p\ra\twoheadrightarrow\la q\ra$ relatively to the cartesian fibration $\F\rightarrow\nabla$.
\end{definition}

\begin{proposition}\label{O_localization}
The inclusion $\con^<(I)\hookrightarrow\F$ admits a left adjoint $\F\rightarrow\con^<(I)$.
\end{proposition}

\begin{lemma}\label{face_fully_faithful}
Let $f:\la p\ra\twoheadrightarrow\la q\ra$ be a face map. Then the induced functor $\enter{\pos{q}}{\Delta^q}\rightarrow\enter{\pos{p}}{\Delta^p}$ is fully faithful.
\end{lemma}

\begin{proof}
It follows from \cite[Lemma 2.10]{volpe} upon verifying that the map $\pos{q}\rightarrow \pos{p}$  is injective and full, and that the square
\[
\begin{tikzcd}
\Delta^q\arrow[d]\arrow[r,"f^*"] & \Delta^p\arrow[d] \\
\pos{q}\arrow[r,"f^*"'] & \pos{p}
\end{tikzcd}
\]
is a pullback square of topological spaces.

For the first claim, let $(g:\la q\ra\twoheadrightarrow\la r\ra),(g':\la q\ra\twoheadrightarrow\la r'\ra)\in \pos{q}$ such that $f^*g=f^*g'\in\pos{p}$. By definition, we have $g\circ f=g'\circ f$ whence $g=g'$ because $f$ is surjective. For the fullness, let $(g:\la q\ra\twoheadrightarrow\la r\ra),(g':\la q\ra\twoheadrightarrow\la r'\ra)\in \pos{q}$ such that $f^*g\leq f^*g'$.  Then there exists a face map $\theta:\la r'\ra\twoheadrightarrow\la r\ra$ such that $g\circ f=\theta\circ g'\circ f$. Since $f$ is surjective, one deduces that $g=\theta\circ g'$ whence $g\leq g'$.  

For the second claim, observe that the map $\Delta^q\rightarrow \Delta^p$ is an embedding whose image consists of all points $x:\la p\ra\rightarrow I$ factoring through the face map $f:\la p\ra\twoheadrightarrow\la q\ra$. By inspection, the image is the fiber product $\pos{q}\times_{\pos{p}}\Delta^p$ and thus the square is indeed a pullback square.
\end{proof}

\begin{lemma}
The objects of $\con^<(I)$ are $\cart(\nabla^{surj})$-local.
\end{lemma}

\begin{proof}
Let $(u:\la q\ra\hookrightarrow I)\in \con^<(I)$ and $\alpha:x\rightarrow y$ a cartesian lift of a face map $f:\la p'\ra\twoheadrightarrow\la p\ra$. Without loss of generality, we identify $x$ with $f^*y$ and the induced path $f^*y\rightarrow f^*y$ in $\exit{\pos{p'}}{\Delta^{p'}}$ is the constant one. Then the induced map $\F(y,u)\rightarrow\F(f^*y,u)$ is weakly homotopy equivalent to the map
\[
\bigsqcup\limits_{g:\la p\ra\rightarrow\la q\ra}\text{Enter}^{\pos{p}}(\Delta^p)(y,g^*u)\xrightarrow{f^*}\bigsqcup\limits_{g':\la p'\ra\rightarrow\la q\ra}\text{Enter}^{\pos{p'}}(\Delta^{p'})(f^*y,(g')^*u)
\]
and it is enough to prove the homotopy equivalence componentwise. Fix a map $g':\la p'\ra\rightarrow\la q\ra$. We distinguish two cases:
\begin{enumerate}
\item $g'$ factors through $f$: write $g'=gf$ with $g:\la p\ra\rightarrow\la q\ra$. Note that $g$ is unique because $f$ is surjective. \Cref{face_fully_faithful} asserts that the functor $f^*:\enter{\pos{p}}{\Delta^p}\rightarrow\enter{\pos{p'}}{\Delta^{p'}}$ is fully faithful so that we get a weak homotopy equivalence 

\[
\text{Enter}^{\pos{p}}(\Delta^p)(y,g^*u)\rightarrow\text{Enter}^{\pos{p'}}(\Delta^{p'})(f^*y,f^*g^*u).
\]

\item $g'$ does not factor through $f$: let $\theta\in \pos{p'}$ the stratum where $f^*y$ lies. By definition, the map $f^*y:\la p'\ra\rightarrow I$ factors through $f$ which implies $\theta\leq f$. Now let $\theta'$ be the stratum where $(g')^*u$ lies. Because $g'$ does not factor through $f$ and $u$ is injective, we have $\theta'\not\leq f$. In particular $\theta'\not\leq\theta$ and therefore $\enter{\pos{p'}}{\Delta^{p'}}(f^*y,(g')^*u)=\emptyset$.
\end{enumerate}

To summarize, we have proved that each component on the right side is either empty, or is homotopy equivalent to exactly one term on the left side. Therefore the map is indeed a homotopy equivalence.
\end{proof}

\begin{proof}[Proof of \Cref{O_localization}]
By \cite[Proposition 5.1.10]{land}, it is enough to specify for each object $(t:\la p\ra\rightarrow I)\in\F$, an object $L(t)\in \con^<(I)$ and a morphism $\alpha_t: t\rightarrow\iota(L(t))$ for which the composition

\begin{equation}\tag{\ensuremath{\ast}}\label{unit_localization}
\con^<(I)(L(t),u)\xrightarrow{\iota}\F(\iota(L(t)),\iota(u))\xrightarrow{\alpha_t^*}\F(t,\iota(u))
\end{equation}

is a weak homotopy equivalence for any object $u\in\con^<(I)$. Consider the factorization
\[
\begin{tikzcd}
\la p\ra\arrow[dr,"t"']\arrow[rr,"f",two heads] & & \la p'\ra\arrow[dl,"L(t)",hook] \\
& I & 
\end{tikzcd}
\]
by a surjection followed by an injection. Define $\alpha_t:t\rightarrow L(t)$ to be the cartesian lift of $f$ with codomain $L(t)$ relatively to the cartesian fibration $\F\rightarrow\nabla$. Then $\alpha_t\in\cart(\nabla^{surj})$ and $L(t)\in \con^<(I)$.

Thus the map $\con^<(I)(L(t),u)\rightarrow\F(\iota(L(t)),\iota(u))$ is a weak homotopy equivalence because the inclusion $\iota:\con^<(I)\hookrightarrow\F$ is fully faithful. Moreover, the map $\F(\iota(L(t)),\iota(u))\xrightarrow{\alpha_t^*}\F(t,\iota(u))$ is also a weak homotopy equivalence because $L(t)$ is $\cart(\nabla^{surj})$-local. Therefore the map \eqref{unit_localization} is also a weak homotopy equivalence.
\end{proof}

\begin{proposition}\label{O_equiv_nabla}
The functor $\con^<(I)\hookrightarrow\F\rightarrow\nabla$ is an equivalence of $\infty$-categories.
\end{proposition}

\begin{proof}
The functor is essentially surjective, hence it remains to prove it is fully faithful. Let $u:\la p\ra\hookrightarrow I,v:\la q\ra\hookrightarrow I\in\con^<(I)$ and $\theta:\la p\ra\rightarrow\la q\ra$ a morphism in $\nabla$. Because $\con^<(I)$ is a full subcategory of $\F$ and $\F\rightarrow\nabla$ is a cartesian fibration, the fiber of the map $\con^<(I)(u,v)\rightarrow\nabla(\la p\ra,\la q\ra)$ over $\theta$ is equivalent to $\enter{\pos{p}}{\Delta^p}(u,\theta^*v)$. By \cite[Theorem A.6.10]{ha}, this mapping space is empty or weakly contractible.

Because $u$ is in the top stratum of $(\Delta^p,\pos{p})$ the fiber is nonempty hence it is weakly contractible. We conclude that the map $\con^<(I)(u,v)\rightarrow\nabla(\la p\ra,\la q\ra)$ is a weak homotopy equivalence.
\end{proof}

The upshot of \Cref{O_equiv_nabla} is that a simplicial object in $\con(M)$ is equivalent to a functor $\con^<(I)\rightarrow\con(M)$.

\section{The functor $\gamma:\con^<(I)\rightarrow \con(M)$}\label{section_gamma}
In this section, we fix a finite simplicial set $X$ whose realization is homeomorphic to $ M$. 

Moreover for any object $\us{n}\in\Fin$, we get a new simplicial set $X^n\cong\Fin(\us{n},X):\nabla\rightarrow\Fin$. This gives a functor $\Fin^{op}\rightarrow\Fun(\nabla,\Fin),\us{n}\mapsto X^n$. Each $X^n$ is a finite simplicial set so that the stratification $\pi_{X^n}:\norm{X^n}\rightarrow\pos{X^n}$ is conical by \cite[Proposition 4.1]{stratified_haya}. This provides a well-defined functor $\enter{\pos{X^-}}{\norm{X^-}}:\Fin^{op}\rightarrow\Cat$ sending an object $\us{n}$ to $\enter{\pos{X^n}}{\norm{X^n}}$.

We thus get a (noncommutative) diagram 

\[
\begin{tikzcd}
\nabla^{op}\arrow[rr,"\text{Enter}^{\pos{\bullet}}{(\Delta^\bullet})"]\arrow[d,"X"'] & & \Cat \\
\Fin^{op}\arrow[urr,"\enter{\pos{X^-}}{\norm{X^-}}"']
\end{tikzcd}.
\]

The goal of this section is to define a natural transformation 
\[
\psi:\text{Enter}^{\pos{\bullet}}(\Delta^\bullet)\Rightarrow \text{Enter}^{\pos{-}}(\norm{X^-})\circ X.
\]


\begin{definition}
For each $p\geq 0$, the stratified map 
\[
\eta_p:(\pi_p:\Delta^p\rightarrow \pos{p})\rightarrow (\pi_{X^p}:\norm{X^{X_p}}\rightarrow\pos{X^{X_p}})
\]
is defined by the composition
\[
(\Delta^p\rightarrow \pos{p})\hookrightarrow\bigsqcup\limits_{\sigma\in\Delta(X^{X_p})}(\Delta^{\sigma}\rightarrow \pos{\sigma})\rightarrow (\norm{X^{X_p}}\rightarrow\pos{X^{X_p}})
\]
where the inclusion $(\Delta^p\rightarrow \pos{p})\hookrightarrow\bigsqcup\limits_{\sigma\in\Delta(X^{X_p})}(\Delta^{\sigma}\rightarrow \pos{\sigma})$ sends the source to the component $(X_p=X_p)\in X^{X_p}_p$.
\end{definition}

\begin{lemma}
The maps $\{\eta_p\}_{p\geq 0}$ extend to a natural transformation $\eta:(\pi_\bullet:\Delta^\bullet\rightarrow\pos{\bullet})\Rightarrow (\pi_{X^-}:\norm{X^-}\rightarrow\pos{X^-})\circ X$ as functors $\nabla^{op}\rightarrow\Strat$.

\end{lemma}

\begin{proof}
Because the second map of the definition is natural in $p$, it is enough to prove that the first one as well. Given a morphism $\theta:\la q\ra\rightarrow\la p\ra$, consider the square 
\[
\begin{tikzcd}
\Delta^p\arrow[d,"\theta^*"']\arrow[r] & \Delta^{(X_p=X_p)}\arrow[d] \\
\Delta^q\arrow[r] & \Delta^{(X_q=X_q)}
\end{tikzcd}.
\]
The right vertical map is given by the composition $\Delta^{(X_p=X_p)}\rightarrow \Delta^{(X_q\rightarrow X_p)}\xrightarrow{\theta^*}\Delta^{(X_q=X_q)}$ with the first map given by the morphism $X^{X_p}\xrightarrow{\theta}X^{X_q}$. 

Thus both the upper right and lower left compositions map $(t:\la p\ra\rightarrow I)$ to $(t\theta:\la q\ra\rightarrow\la p\ra\rightarrow I)\in\Delta^{(X_q=X_q)}$.
\end{proof}

\begin{definition}
The natural transformation $\psi:\enter{\pos{\bullet}}{\Delta^\bullet}\Rightarrow \enter{\pos{-}}{\norm{X^-}}\circ X$ is the composition $\enter{-}{-}\circ\eta$.
\end{definition}

\begin{definition}
The functors 
\begin{align*}
\mathcal{Q}\rightarrow\nabla, \\
p:\B_X\rightarrow\Fin
\end{align*}
are the respective weighted nerves of the functors 
\begin{align*}
\text{Enter}^{\pos{X^-}}(\norm{X^-})\circ X:\nabla^{op}\rightarrow\Cat^\delta, \\
\text{Enter}^{\pos{X^-}}(\norm{X^-}):\Fin^{op}\rightarrow\Cat.
\end{align*}
\end{definition}

We get a commutative diagram

\[
\begin{tikzcd}
\F\arrow[dr]\arrow[r] & \mathcal{Q}\arrow[d]\arrow[r]\arrow[dr,near start, phantom,"\lrcorner"] & \B_X\arrow[d,"p"]\\
& \nabla\arrow[r,"X"'] & \Fin 
\end{tikzcd}
\]

where the functor $\F\rightarrow\mathcal{Q}$ is induced by the natural transformation $\psi$.

\begin{definition}
We define $\delta:\F\rightarrow\B$ to be the composition $\F\rightarrow\mathcal{Q}\rightarrow\B_X$.
\end{definition}

Moreover, \Cref{refinement_of_partitions} provides a natural transformation $\enter{\pos{X^-}}{\norm{X^-}}\Rightarrow\enter{\prt{\us{-}}}{M^-}$. This gives a cartesian functor
\[
\begin{tikzcd}
\B_X\arrow[dr,"p"']\arrow[rr,"\mu"] & & \X_M \arrow[dl] \\
& \Fin & 
\end{tikzcd}
\]

and composing with $\delta$ provides a functor $\bar{\gamma}:\F\rightarrow\X_M$.

\begin{lemma}
The functor $\bar{\gamma}:\F\rightarrow\X_M$ restricts to a functor $\gamma:\con^<(I)\rightarrow\con(M)$.
\end{lemma}

\begin{proof}
Let $u:\la p\ra\hookrightarrow I$ an object of $\con^<(I)$. By definition, the object $\bar{\gamma}(u)$ is identified with the image of $u$ under the composition $\Delta^p\xrightarrow{\eta_p}\norm{X^{X_p}}\xrightarrow{\cong}M^{X_p}$. This composition extends to a commutative rectangle
\[
\begin{tikzcd}
\Delta^{p}\arrow[r,"\eta_p"]\arrow[d] & \norm{X^{X_p}}\arrow[d]\arrow[r] & M^{X_p}\arrow[d] \\
\pos{p}\arrow[r] & \pos{X^{X_p}}\arrow[r] & \prt{X_p}
\end{tikzcd}.
\]
By definition, $u$ is sent to the top stratum $1_{\la p\ra}=(\la p\ra = \la p\ra)\in \pos{p}$. Then $\eta_p(1_{\la p\ra})=(X_p=X_p)\in X^{X_p}_p$ by definition of $\eta_p$. Because the map $X_p=X_p$ is injective, this element is in turn mapped to the top element in $\prt{X_p}$. 

To summarize, the element $u\in\Delta^p$ is mapped to a point lying in the top stratum, hence is a configuration.
\end{proof}

\section{Cofinality of the functor $\gamma$}\label{section_cofinality}

Let $\us{n}\in\Fin$ and consider the diagram

\[
\begin{tikzcd}[column sep=2.6em, row sep=2.6em]
(X_{\us{n}/})^{op}
  \arrow[dd]
  \arrow[drrrr,bend left=15,"\enter{\pos{\bullet}}{\Delta^\bullet}"]
  \arrow[dr]
& & & & \\
&
\nabla^{op}
  \arrow[rrr,"\enter{\pos{\bullet}}{\Delta^\bullet}"]
  \arrow[dd,"X"',crossing over]
&
{} \arrow[dd, phantom,
  "\rotatebox{-90}{\(\Longrightarrow\)}"{description,pos=.35},
  "\psi"{right,pos=.35}
]
& {} &
\Cat^\delta \\
(\Fin_{\us{n}/})^{op}
  \arrow[dr]
  \arrow[urrrr]
& & & & \\
&
\Fin^{op}
  \arrow[uurrr,"\enter{\pos{X^-}}{\norm{X^-}}"']
& {} & &
\end{tikzcd}.
\]

The natural transformation $\psi:\enter{\pos{\bullet}}{\Delta^\bullet}\Rightarrow\enter{\pos{X^-}}{\norm{X^-}}\circ X$ in $\Fun(\nabla^{op},\Cat^\delta)$ defines, by precomposition, a natural transformation $\psi_n:\enter{\pos{\bullet}}{\Delta^\bullet}\Rightarrow\enter{\pos{X^-}}{\norm{X^-}}\circ X$ in $\Fun((X_{\us{n}/})^{op},\Cat^\delta)$.

Now observe that the functor $\Fin_{\us{n}/}\rightarrow\Fin$ can be identified as the unstraightening of the functor $\Fin(\us{n},-):\Fin\rightarrow\Set\hookrightarrow\cS$. Precomposition with the functor $X:\nabla\rightarrow\Fin$ gives precisely $X^n\cong(\Fin(\us{n},-)\circ X):\nabla\rightarrow\Set\hookrightarrow\cS$.
The unstraightening of this functor is then equivalent to the projection  $\Delta(X^n)^{op}\rightarrow\nabla$. We then deduce the equivalence $\Delta(X^n)^{op}\simeq X_{\us{n}/}$. 

\begin{definition}
For $n\geq 0$, we define 
\begin{align*}
\cE_n\rightarrow\Delta(X^n)^{op} \\
P_n\rightarrow \Fin_{\us{n}/}
\end{align*} 
to be the respective weighted nerves of the functors 
\begin{align*}
\text{Enter}^{\pos{\bullet}}(\Delta^\bullet):\Delta(X^n)\rightarrow\Cat^\delta, \\
\text{Enter}^{\pos{X^-}}(\norm{X^-}):(\Fin_{\us{n}/})^{op}\rightarrow\Cat^\delta.
\end{align*}
\end{definition}

The weighted nerves $\cE_n,P_n$ fit into a commutative diagram

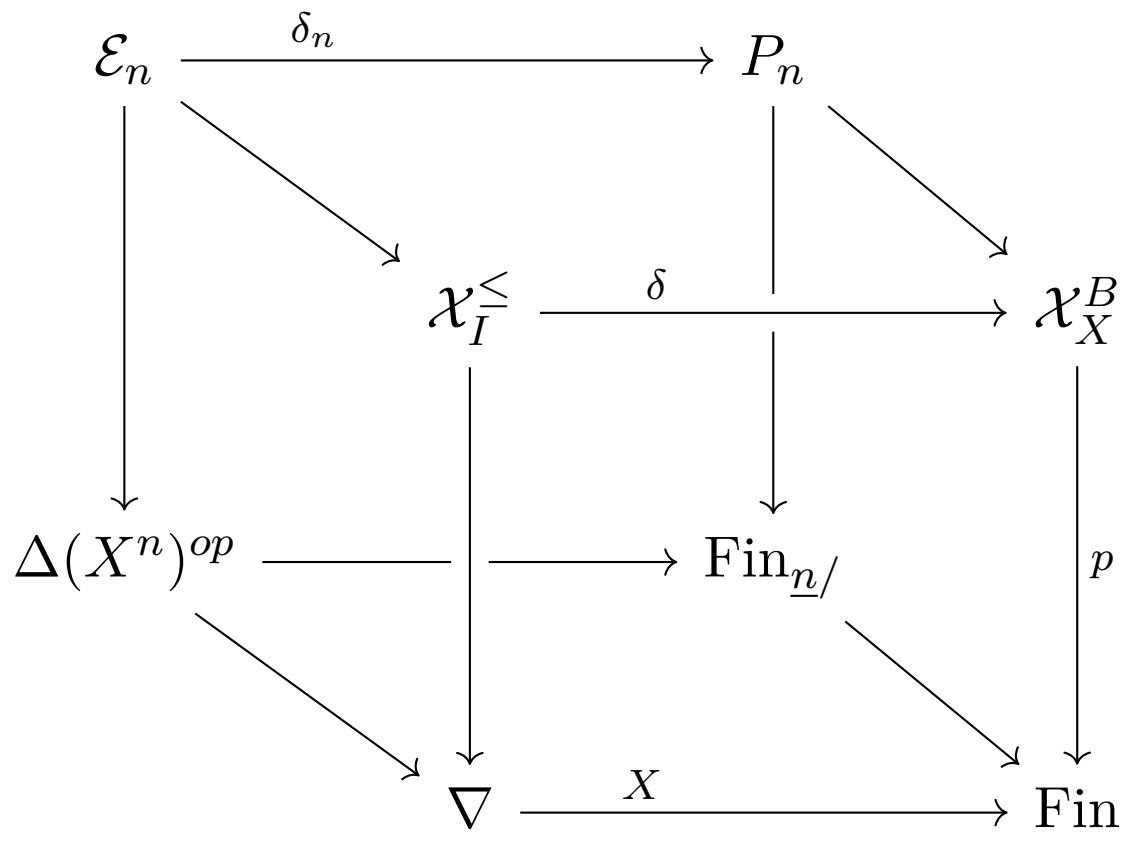
\[
\begin{tikzcd}[row sep=2.5em]
\cE_n \arrow[rr,"\delta_n",near start] \arrow[dr] \arrow[dd] &&
  P_n \arrow[dd] \arrow[dr] \\
& \F \arrow[rr,crossing over,"\delta",near start] &&
  \B_X \arrow[dd,"p"] \\
\Delta(X^n)^{op} \arrow[rr] \arrow[dr] && \Fin_{\us{n}/} \arrow[dr] \\
& \nabla \arrow[rr,"X",near start] \arrow[uu,<-,crossing over]&& \Fin
\end{tikzcd}
\]

where all vertical arrows are cartesian fibrations, the left, lower and right faces are pullbacks and the functor $\delta_n:\cE_n\rightarrow P_n$ is induced from the natural transformation $\psi_n$. Then, one gets that the top face is also a pullback.

Next consider the diagram where all squares are pullbacks:

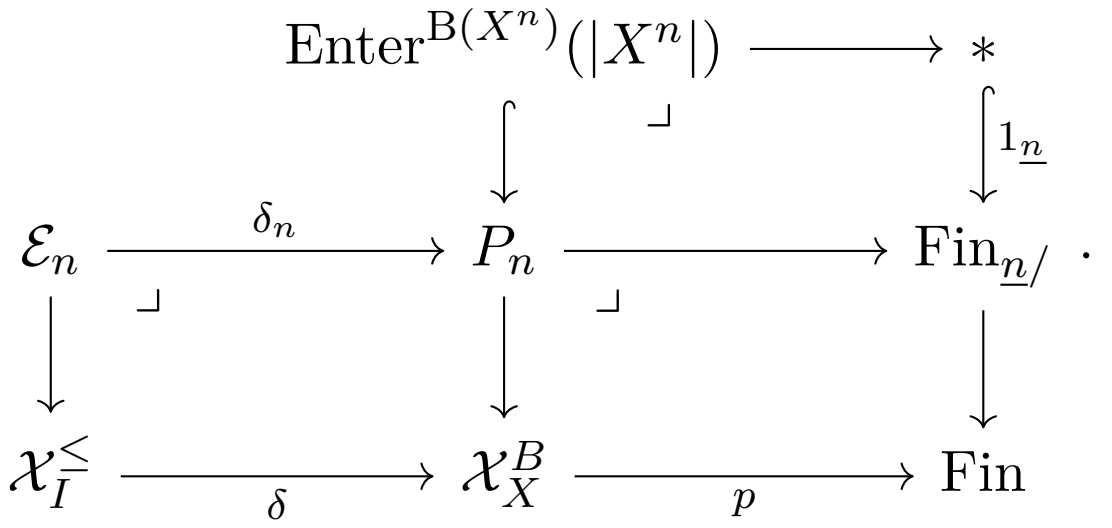
\[
\begin{tikzcd}
& \enter{\pos{X^n}}{\norm{X^n}}\arrow[d,hook]\arrow[r]\arrow[dr, phantom, "\lrcorner", very near start] & \ast\arrow[d,"1_{\us{n}}",hook] \\
\cE_n\arrow[d]\arrow[r,"\delta_n"]\arrow[dr, phantom, "\lrcorner", very near start] & P_n\arrow[d]\arrow[r]\arrow[dr, phantom, "\lrcorner", very near start] & \Fin_{\us{n}/}\arrow[d] \\
\F\arrow[r,"\delta"'] & \B_X\arrow[r,"p"'] & \Fin
\end{tikzcd}.
\]

The $\infty$-category $\Fin_{\us{n}/}$ admits an initial object, the identity $1_{\us{n}}$. By \Cref{adjoint_cart}, the fiber inclusion $\enter{\pos{X^n}}{\norm{X^n}}\simeq p^{-1}({\us{n}})\hookrightarrow P_n$ has a right adjoint we call $\rho_n$. Define $\lambda_n:\cE_n\rightarrow \enter{\pos{X^n}}{\norm{X^n}}$ to be the composition $\cE_n\xrightarrow{\delta_n} P_n\xrightarrow{\rho_n} \enter{\pos{X^n}}{\norm{X^n}}$. Finally, fix an object $x\in\B_X$ such that $p(x)=\us{n}$ and consider the pullbacks

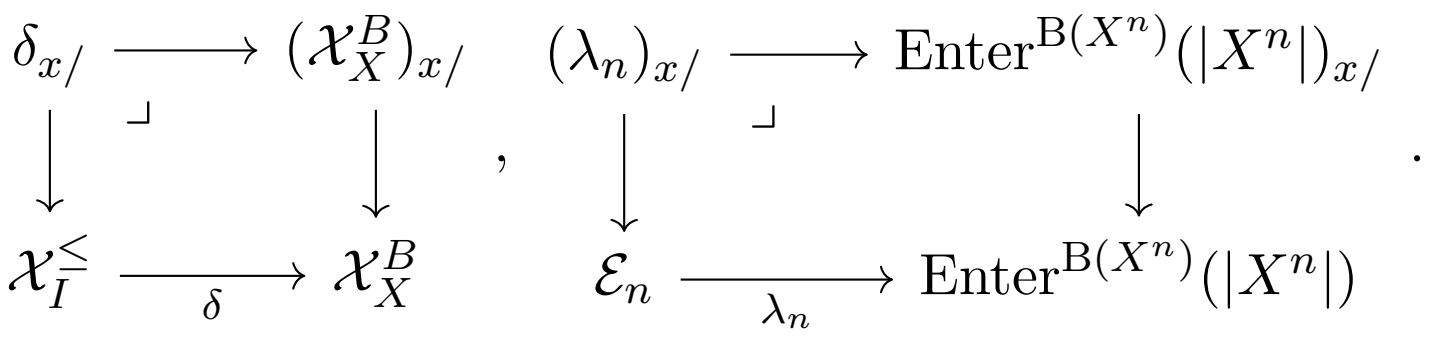
\[
\begin{tikzcd}
\delta_{x/}\arrow[r]\arrow[d]\arrow[dr, phantom, "\lrcorner", very near start] & (\B_X)_{x/}\arrow[d] \\
\F\arrow[r,"\delta"'] & \B_X
\end{tikzcd},
\begin{tikzcd}
(\lambda_n)_{x/}\arrow[r]\arrow[d]\arrow[dr, phantom, "\lrcorner", very near start] & \enter{\pos{X^n}}{\norm{X^n}}_{x/}\arrow[d] \\
\cE_n\arrow[r,"\lambda_n"'] & \enter{\pos{X^n}}{\norm{X^n}}
\end{tikzcd}.
\]

Notice that the right square naturally maps to the left square so that we get a natural functor $(\lambda_n)_{x/}\rightarrow\delta_{x/}$.

\begin{lemma}\label{equivalence_slices}
The functor $(P_n)_{x/}\rightarrow(\B_X)_{x/}$ is an equivalence of $\infty$-categories.
\end{lemma}

\begin{proof}
By definition of $P_n=\B_X\times_{\Fin}\Fin_{\us{n}/}$, we get a pullback of slices 

\[
\begin{tikzcd}
(P_n)_{x/}\arrow[d]\arrow[r] & (\Fin_{\us{n}/})_{1_{\us{n}}/}\arrow[d] \\
(\B_X)_{x/}\arrow[r] & \Fin_{\us{n}/}
\end{tikzcd},
\]

hence it is enough to prove that we have an equivalence $(\Fin_{\us{n}/})_{1_{\us{n}}/}\rightarrow \Fin_{\us{n}/}$. Because this map is a left fibration, it is enough to check that all fibers are contractible spaces. But the fiber over an object $f\in \Fin_{\us{n}/}$ is the mapping space $\Fin_{\us{n}/}(1_{\us{n}},f)$ which is contractible because $1_{\us{n}}$ is an initial object.

\end{proof}

\begin{proposition}\label{equivalence_slices_pullback}
The functor $(\lambda_n)_{x/}\rightarrow\delta_{x/}$ is an equivalence of $\infty$-categories.
\end{proposition}

\begin{proof}
Consider first the rectangle 

\[
\begin{tikzcd}
(\delta_n)_{x/}\arrow[r]\arrow[d] & (P_n)_{x/}\arrow[r]\arrow[d] & \enter{\pos{X^n}}{\norm{X^n}}_{x/}\arrow[d] \\
\cE_n\arrow[r,"\delta_n"'] & P_n \arrow[r,"\rho_n"'] & \enter{\pos{X^n}}{\norm{X^n}}
\end{tikzcd}.
\]

By \Cref{adjoint_cart} and \Cref{pullback_adj}, the right square is a pullback. The left square is also a pullback by definition so that the whole rectangle is a pullback. This gives an equivalence $(\delta_n)_{x/}\simeq(\lambda_n)_{x/}$.

Now consider the cube

\[
\begin{tikzcd}[row sep=2.5em]
(\delta_n)_{x/} \arrow[rr] \arrow[dr] \arrow[dd] &&
  (P_n)_{x/} \arrow[dd] \arrow[dr] \\
& \cE_n \arrow[rr,crossing over,"\delta_n",near start] &&
  P_n \arrow[dd] \\
\delta_{x/} \arrow[rr] \arrow[dr] && (\B_X)_{x/} \arrow[dr] \\
& \F \arrow[rr,"\delta",near start] \arrow[uu,<-,crossing over]&& \B_X
\end{tikzcd}.
\]

The upper, front and lower faces are pullbacks and thus so is the back one. By \Cref{equivalence_slices}, the functor $(P_n)_{x/}\rightarrow(\B_X)_{x/}$ is an equivalence of categories and thus so is its pullback $(\delta_n)_{x/}\rightarrow\gamma_{x/}$.

Composing the equivalences, one gets
\[
(\lambda_n)_{x/}\simeq(\delta_n)_{x/}\xrightarrow{\simeq}\delta_{x/}.
\]

\end{proof}

\begin{corollary}\label{balancing_lemma}
The functor $\delta:\F\rightarrow\B_X$ is cofinal if and only if $\lambda_n:\cE_n\rightarrow \enter{\pos{X^n}}{\norm{X^n}}$ is cofinal for all integers $n\geq 0$.
\end{corollary}

\begin{proof}
By Quillen's theorem A, both cofinality are equivalent to the weak  contractibility of the pullbacks $\delta_{x/},(\lambda_n)_{x/}$ for all objects $x\in\B_X$. These pullbacks are equivalent by \Cref{equivalence_slices_pullback}.
\end{proof}

Next we want to better understand the functor 
\[
\lambda_n:\cE_n\rightarrow\enter{\pos{X^n}}{\norm{X^n}}.
\]

\begin{lemma}\label{lambda_refraction}
The functor $\lambda_n:\cE_n\rightarrow\enter{\pos{X^n}}{\norm{X^n}}$ is equivalent to the refraction diagram of the cocone $\Delta(X^n)^\rhd\rightarrow\Cat^\delta,\sigma\mapsto\enter{\pos{\sigma}}{\Delta^\sigma},\infty\mapsto\enter{\pos{X^n}}{\norm{X^n}}$.
\end{lemma}

\begin{proof}
It is enough to prove that for each simplex $\sigma\in\Delta(X^n)$, the functor $\enter{\pos{\sigma}}{\Delta^\sigma}\rightarrow\enter{\pos{X^n}}{\norm{X^n}}$ obtained by contravariant transport is equivalent to the functor $\enter{\pos{\sigma}}{\Delta^\sigma}\rightarrow\enter{\pos{X^n}}{\norm{X^n}}$ induced by the successive maps
\[
(\Delta^\sigma\rightarrow\pos{\sigma})\hookrightarrow\bigsqcup\limits_{\tau:\us{n}\rightarrow X}(\Delta^\tau\rightarrow\pos{\tau})\rightarrow(\norm{X^n}\rightarrow\pos{X^n}).
\]

By definition of $\psi_n$, the functor $\gamma_n$ is induced for each $p$-simplex $\sigma\in X^n_p$ by the composition 

\[
(\Delta^\sigma\rightarrow\pos{\sigma})\hookrightarrow\bigsqcup\limits_{\tau:X_p\rightarrow X}(\Delta^\tau\rightarrow\pos{\tau})\rightarrow\lf(\norm{X^{X_p}}\rightarrow\pos{X^{X_p}}\ri).
\]

Moreover, because the functor $\rho_n:P_n\rightarrow\enter{\pos{X^n}}{\norm{X^n}}$ is a refraction diagram, the induced functor on the fiber over $\sigma\in\Delta(X^n)$ is induced by the map

\[
\lf(\norm{X^{X_p}}\rightarrow\pos{X^{X_p}}\ri)\xrightarrow{\sigma:\us{n}\rightarrow X_p}(\norm{X^{n}}\rightarrow\pos{X^{n}}).
\]

Since $\lambda_n=\rho_n\circ\gamma_n$, one deduces that the induced functor $\enter{\pos{\sigma}}{\Delta^\sigma}\rightarrow\enter{\pos{X^n}}{\norm{X^n}}$ is induced by the composition

\begin{align*}
(\Delta^\sigma\rightarrow\pos{\sigma}) \hookrightarrow\bigsqcup\limits_{\tau:X_p\rightarrow X}(\Delta^\tau\rightarrow\pos{\tau})
\rightarrow(\norm{X^{X_p}}\rightarrow\pos{X^{X_p}})
\xrightarrow{\sigma:\us{n}\rightarrow X_p}(\norm{X^{n}}\rightarrow\pos{X^{n}})
\end{align*}

which is precisely the composition

\[
(\Delta^\sigma\rightarrow\pos{\sigma})\hookrightarrow\bigsqcup\limits_{\tau:\us{n}\rightarrow X}(\Delta^\tau\rightarrow\pos{\tau})\rightarrow(\norm{X^n}\rightarrow\pos{X^n}).
\]

\end{proof}

\begin{corollary}\label{lambda_cofinal}
For each $n\geq 0$, the functor 
\[
\lambda_n:\cE_n\rightarrow\enter{\pos{X^n}}{\norm{X^n}}
\]
is cofinal.
\end{corollary}

\begin{proof}
By \Cref{lambda_refraction}, $\lambda_n$ is the refraction diagram of the cocone $\Delta(X^n)^\rhd\rightarrow\Cat$. By \cite[Theorem 1.2]{stratified_haya}, this is a colimit diagram. By \cite[Theorem 7.4.5.13]{kerodon}, we deduce that $\lambda_n$ is a localization functor, hence is cofinal.
\end{proof}

To conclude this section, consider the diagram

\[
\begin{tikzcd}
\con^<(I)\arrow[d,hook]\arrow[rr,"\gamma"] & & \con(M)\arrow[d,hook] \\
\F\arrow[d]\arrow[r,"\delta"] & \B_X\arrow[r,"\mu"]\arrow[dr] & \X_M\arrow[d] \\
\nabla\arrow[rr,"X"'] & & \Fin  
\end{tikzcd}.
\]

\begin{lemma}\label{mu_cofinal}
The functor $\mu:\B_X\rightarrow\X_M$ is cofinal.
\end{lemma}

\begin{proof}
The functor $\mu$ is cartesian and induces on fibers over $\us{n}\in\Fin$ the functor $\enter{\pos{X^n}}{\norm{X^n}}\rightarrow\enter{\prt{\us{n}}}{M^n}$. By \cite[Theorem 0.3.1]{exodromy}, this functor is a localization functor. Because $\mu$ induces localizations on fibers, \cite[\href{https://kerodon.net/tag/02LW}{Proposition 6.3.4.2}]{kerodon} tells us that $\mu$ itself is a localization functor, hence is cofinal.
\end{proof}

\begin{lemma}\label{delta_cofinal}
The functor $\delta:\F\rightarrow\B_X$ is cofinal.
\end{lemma}

\begin{proof}
By \Cref{balancing_lemma}, it is enough to prove that the functor $\lambda_n:\cE_n\rightarrow\enter{\pos{X^n}}{\norm{X^n}}$ is cofinal for all $n\geq 0$. This is the case by \Cref{lambda_cofinal}.
\end{proof}

\begin{proposition}\label{gamma_cofinal}
The functor $\gamma:\con^<(I)\rightarrow\con(M)$ is cofinal.
\end{proposition}

\begin{proof}
By \Cref{cofinal_ff}, it is enough to prove that the composition 
\[
\con^<(I)\xrightarrow{\gamma}\con(M)\hookrightarrow\X_M
\]
is cofinal. This is nothing but the composition 
\[
\con^<(I)\hookrightarrow\F\xrightarrow{\delta}\B\xrightarrow{\mu}\X_M.
\] 
The first functor is a right adjoint by \Cref{O_localization} hence is cofinal. The second and third ones are cofinal by \Cref{mu_cofinal} and \Cref{delta_cofinal}. Therefore their composition is also cofinal.

\end{proof}

\section{The main result}
In this short section, we assemble all the relevant results to prove the main result.

\begin{theorem}(\Cref{main_simplicial_model})
Let $M$ be a closed $d$-dimensional smooth manifold and $X$ a finite simplicial set whose geometric realization is homeomorphic to $M$. Then there exists a cofinal functor $\Gamma_X:\Delta^{op}\rightarrow\cDisk_{/M}$ making the following triangle commute:
\[
\begin{tikzcd}
\Delta^{op}\arrow[dr,"X"']\arrow[r,"\Gamma_X"] & \cDisk_{/M}\arrow[d,"\pi_0"] \\
& \Fin
\end{tikzcd}.
\]
\end{theorem}

\begin{proof}
Consider the diagram 

\[
\begin{tikzcd}
\con^<(I)\arrow[d]\arrow[r,"\gamma"] & \con(M)\arrow[d] \\
\nabla\arrow[r,"X"'] & \Fin
\end{tikzcd}.
\]

By \Cref{gamma_cofinal}, the functor $\gamma:\con^<(I)\rightarrow\con(M)$ is cofinal. Moreover by \Cref{O_equiv_nabla}, the functor $\con^<(I)\rightarrow\nabla\cong\Delta^{op}$ is an equivalence of $\infty$-categories and by \Cref{equiv_disks}, we have another equivalence $\con(M)\simeq\cDisk_{/M}$ in $(\Cat)_{/\Fin}$. Thus we get a cofinal functor $\Gamma_X:\Delta^{op}\rightarrow\cDisk_{/M}$ fitting into a commutative triangle 

\[
\begin{tikzcd}
\Delta^{op}\arrow[dr,"X"']\arrow[r,"\Gamma_X"] & \cDisk_{/M}\arrow[d,"\pi_0"] \\
& \Fin
\end{tikzcd}.
\]

\end{proof}

To conclude this section, we observe that \Cref{main_simplicial_model} generalizes immediately to manifolds with tangential structures. Following \cite[Section 2]{factorization}, for a topological group $G$ and a map $\theta:BG\rightarrow BO(d)$, a $\theta$-structure on a smooth $d$-dimensional manifold $M$ without boundary is a lifting, up to homotopy, of its tangent classifier $\tau_M:M\rightarrow BO(d)$

\[
\begin{tikzcd}
& BG\arrow[d,"\theta"] \\
M\arrow[ur,"\tau_M^g"]\arrow[r,"\tau_M"'] & BO(d)
\end{tikzcd}.
\]

By abuse of notation, we will leave the map $\theta:BG\rightarrow BO(d)$ implicit and talk about $G$-structures. Manifolds with $G$-structures with embeddings between them preserving the $G$-structure assemble into an $\infty$-category $\cMan^G$ formally defined (\cite[Definition 2.7]{factorization}) by the pullback of $\infty$-categories
\[
\begin{tikzcd}
\cMan^G\arrow[d]\arrow[r]\arrow[dr,"\lrcorner",phantom,very near start] & \cS_{/BG}\arrow[d,"\theta"] \\
\cMan\arrow[r,"\tau_{-}"'] & \cS_{/BO(d)}
\end{tikzcd}.
\]

For objects $M,N\in\cMan^G$, we write $\Emb^G(M,N):=\cMan^G(M,N)$.
We get the full subcategory $\cDisk^G\subset\cMan^G$ and for a manifold $M\in\cMan^G$, the overcategory $\cDisk^G_{/M}$. 

The remarkable observation in \cite[Lemma 2.5]{factorization} is that the inclusion $\cDisk^G_{/M}\rightarrow\cDisk_{/M}$ is an equivalence of $\infty$-categories. Because of this, we have a straightforward generalization of \Cref{main_simplicial_model} in the $G$-structured setting.

\begin{theorem}\label{main_structured}
Let $M$ be a closed $d$-dimensional smooth manifold with a $G$-structure and $X$ a finite simplicial set whose geometric realization is homeomorphic to $M$. Then there exists a cofinal functor $\Gamma_X:\Delta^{op}\rightarrow\cDisk^G_{/M}$ making the following triangle commute:
\[
\begin{tikzcd}
\Delta^{op}\arrow[dr,"X"']\arrow[r,"\Gamma_X"] & \cDisk^G_{/M}\arrow[d,"\pi_0"] \\
& \Fin
\end{tikzcd}.
\]
\end{theorem}

\begin{rem}
While \Cref{main_simplicial_model} and \Cref{main_structured} are stated for smooth manifolds, the proofs adapt verbatim to topological manifolds. Indeed, the definitions of the stratifications $\pi_M^-,\pi_X$, the configuration category $\con(M)$, the functor $\gamma:\con^<(I)\rightarrow\con(M)$ and the proof of the cofinality definitely do not depend on the smooth structure. The only part where the smooth structure could possibly intervene is the proof of the equivalence $\con(M)\simeq\cDisk_{/M}$. However, the proof relies on computing mapping spaces in enter path categories, which does not depend on the smooth structure, and on the results of \cite[Section 2]{factorization} which hold for topological manifolds. In fact, these results are stated for topological manifolds, but also hold for smooth manifolds as observed earlier (\cite[Remark 3.29]{factorization},\cite[Remark A.8]{arakawa}). 
\end{rem}

\section{Applications and further prospects}

\subsection{A cosimplicial model for spaces of embeddings}
This was the original motivation of \cite{model}. Let $M$ be a closed $d$-dimensional smooth manifold with a $G$-structure and $X$ a finite simplicial set whose geometric realization is homeomorphic to $M$. We first show how \Cref{main_structured} provides a simplicial presheaf computing $\Emb^G(-,M)$. Let $\mathcal{Y}:\cDisk^G\hookrightarrow\Psh(\cDisk^G)$ be the Yoneda embedding.

\begin{lemma}\label{emb_colim}
There is an equivalence in $\Psh(\cDisk^G)$
\[
\colim{U\in\cDisk^G_{/M}}{\Emb^G(-,U)}\simeq\Emb^G(-,M).
\]
\end{lemma}

\begin{proof}
Consider the composition
\[
(\cDisk^G)^{op}\xrightarrow{\mathcal{Y}}\Psh(\cDisk^G)^{op}\xrightarrow{\Psh(\cDisk^G)(\bullet,\Emb^G(-,M))}\cS.
\]
By definition, its unstraightening is equivalent to the right fibration
\[
\mathcal{Y}_{/\Emb^G(-,M)}\rightarrow\cDisk^G.
\]
Moreover, the Yoneda lemma gives a natural equivalence of presheaves $\Psh(\cDisk^G)(\bullet,\Emb^G(-,M))\circ\mathcal{Y}\simeq \Emb^G(-,M)$. Since the unstraightening of this functor is equivalent to the right fibration
\[
\cDisk^G_{/M}\rightarrow\cDisk^G,
\]
one deduces an equivalence of $\infty$-categories $\cDisk^G_{/M}\simeq \mathcal{Y}_{/\Emb^G(-,M)}$ over $\cDisk^G$. Combining this equivalence with the fact that a presheaf is a colimit of its representables, one deduces the equivalences
\[
\colim{U\in\cDisk^G_{/M}}{\Emb^G(-,U)}\simeq\colim{\mathcal{Y}_{/\Emb^G(-,M)}}{\Emb^G(-,U)}\simeq\Emb^G(-,M).
\]
\end{proof} 

Now define the cosimplicial space $\Omega_F$ by the composition
\[
\Delta\xrightarrow{\Gamma_X^{op}}(\cDisk^G_{/M})^{op}\rightarrow(\cDisk^G)^{op}\hookrightarrow(\cMan^G)^{op}\xrightarrow{F}\cS.
\]

This cosimplicial space allows to recover \cite[Theorem 1.1]{model}.
\begin{proposition}(\Cref{cosimplicial_space_from_Gamma})
Given a presheaf $F:(\cMan^G)^{op}\rightarrow\cS$, the cosimplicial space $\Omega_F$ satisfies
\begin{enumerate}
\item $T_\infty F(M) \simeq \limi{\Delta}{\Omega_F}$.
\item For each cosimplicial map $\theta:[p]\rightarrow [q]$, there is an embedding $\iota_\theta:X_q\times\R^d\hookrightarrow X_p\times\R^d$ such that $\pi_0(\iota_\theta)\cong (X_\theta:X_q\rightarrow X_p)$ and a weak equivalence of morphisms:
\[(\Omega_F^\theta:\Omega_F^p\rightarrow\Omega_F^q)\simeq (\iota_\theta^*:F(X_p\times\R^d)\rightarrow F(X_q\times\R^d)). 
\] 
\end{enumerate}
\end{proposition}

\begin{proof}
\Cref{main_structured} provides a cofinal functor $\Gamma_X:\Delta^{op}\rightarrow \cDisk^G_{/M}$, reducing the equivalence of \Cref{emb_colim} to an equivalence 

\[
\Emb^G(-,M)\simeq\colim{\Delta^{op}}{\Emb^G(-,\Gamma_X)}.
\]

Moreover, by naturality of the Yoneda lemma, we have an equivalence of cosimplicial spaces
\[
\Omega_F\simeq\Psh(\cDisk^G)(\Emb^G(-,\Gamma_X),F).
\]

On the one hand, we have the Taylor approximation

\begin{align*}
T_{\infty}F(M) & := \Psh(\cDisk^G)(\Emb^G(-,M),F) \\
&\simeq \Psh(\cDisk^G)(\colim{\Delta^{op}}{\Emb^G(-,\Gamma_X)},F) \\
&\simeq \limi{\Delta}{\Psh(\cDisk^G)(\Emb^G(-,\Gamma_X),F)} \\
&=\limi{\Delta}{\Omega_F}.
\end{align*}

On the other hand, fix a cosimplicial map $\theta[p]\rightarrow\theta[q]$. \Cref{main_structured} provides the appropriate embedding $\iota_\theta:X_q\times\R^d\hookrightarrow X_p\times\R^d$. 
Since $\Gamma_X$ is degreewise representable, the Yoneda lemma yields the weak homotopy equivalence 
\[
(\Omega^\theta:\Omega^p\rightarrow\Omega^q) \simeq (F(\iota_\theta):F(X_p\times\R^d)\rightarrow F(X_q\times\R^d)).
\]
\end{proof}

In particular, let $W$ be a $D$-dimensional manifold without boundary. If $D-d\geq 3$, Goodwillie-Weiss excision estimate asserts that $T_{\infty}(\Emb(-,W))(M)\simeq\Emb(M,W)$. The upshot is that we get a cosimplicial space $\Omega_W$ such that
\begin{enumerate}
\item $\Emb(M,W)\simeq\limi{\Delta}{\Omega_W}$,
\item for each cosimplicial map $\theta:[p]\rightarrow [q]$, there is an embedding $\iota_\theta:X_q\times\R^d\hookrightarrow X_p\times\R^d$ such that $\pi_0(\iota_\theta)\cong (X_\theta:X_q\rightarrow X_p)$ and a weak equivalence of morphisms:
\[(\Omega_W^\theta:\Omega_W^p\rightarrow\Omega_W^q)\simeq (\iota_\theta^*:\Emb(X_p\times\R^d,W)\rightarrow\Emb(X_q\times\R^d,W)). 
\] 
\end{enumerate}

This functor can be made strict through a classical rigidification process, as stated in \cite[Remark 1.2.6.2]{htt} for example. More precisely, because the category $\sSet$ endowed with the Quillen model structure forms a combinatorial simplicial model category whose bifibrant objects are exactly Kan complexes, \cite[Proposition 4.2.4.4]{htt} implies that there exists a \emph{strict} projectively bifibrant functor $\tilde{\Omega}_W:\Delta\rightarrow\sSet$ satisfying the three properties listed above.

\subsection{A simplicial model for factorization homology}
Let $\cC$ be a symmetric monoidal $\infty$-category (for example $\cC=\cS$ with the cartesian product as tensor product). Given a $G$-structured smooth $d$-dimensional manifold $M$ without boundary and a symmetric monoidal functor $A:(\cDisk^G,\sqcup)\rightarrow (\cC,\otimes)$, \emph{factorization homology} of $M$ with coefficients in $A$ is defined (\cite[Definition 3.2]{factorization}) as the colimit
\[
\int_M A\simeq\colim{\cDisk^G_{/M}}{A}.
\]

If $M$ is closed, fixing a finite simplicial set $X$ and a homeomorphism $\norm{X}\cong M$, \Cref{main_structured} provides a cofinal functor $\Gamma_X:\Delta^{op}\rightarrow\cDisk^G_{/M}$ so that we get an equivalence
\[
\int_M A\simeq\colim{\Delta^{op}}{A\circ\Gamma_X}
\]
and for each $p\geq 0$, we have an equivalence
\[
(A\circ\Gamma_X)_p\simeq A(X_p\times\R^d)\simeq A(\R^d)^{\otimes X_p}.
\]

\begin{rem}
The $\infty$-category $\Psh(\cDisk^G)$ admits a symmetric monoidal structure $\otimes$ given by the Day convolution tensor product. Explicitly, for presheaves $F,G:(\cDisk^G)^{op}\rightarrow\cS$ and an object $\us{n}\times\R^d\in\cDisk^G$, we have
\[
(F\otimes G)(\us{n}\times\R^d)\simeq\bigsqcup\limits_{\substack{S,T\subset\us{n}\\ S\sqcup T=\us{n}}} F(S\times\R^d)\times G(T\times\R^d).
\] 
The Yoneda embedding $\mathcal{Y}:(\cDisk^G,\sqcup)\hookrightarrow (\Psh(\cDisk^G),\otimes)$ is then symmetric monoidal (\cite[Corollary 4.8.1.12]{ha},\cite[Remark 4.8.1.13]{ha}). Thus the factorization homology of $\mathcal{Y}$ is well defined.

Using \cite[Theorem 3.24]{factorization} and adapting \cite[Corollary 2.38]{factorization_bis} to $\cDisk^G$, one gets an equivalence of presheaves

\[
\int_M\mathcal{Y}\simeq \Emb^G(-,M)
\]
which recovers \Cref{emb_colim}.
\end{rem}

\subsection{Manifolds with boundary}

We briefly discuss how to adapt \Cref{main_simplicial_model} to manifolds with boundary. We first adapt our definition of $\con(M)$ to the case with boundary, following \cite[Section 6.1]{confcat}. 

\begin{definition}
The (ordinary) category $\Fin_\ast$ has objects finite pointed sets $\us{n}_\ast:=\{\ast,1,\cdots,n\}$ for $n\geq 0$, and pointed maps as morphisms.
\end{definition}

\begin{definition}
For $n\geq 0$, the poset $\prt{\us{n}_\ast}$ has elements equivalence relations $R\subset\us{n}_\ast\times\us{n}_\ast$ and partial order given by reverse inclusion: $R\leq S$ if and only if $S\subset R$.
\end{definition}

Let $M$ be a smooth, compact $d$-dimensional manifold with nonempty boundary $\partial M$. The inclusion $\partial M\hookrightarrow M$ defines a quotient space $M_\ast:M/\partial M$. This comes with a basepoint being the quotient boundary. For $n\geq 0$, one then has a homeomorphism $M_\ast^n\cong\Map_\ast(\us{n}_\ast,M_\ast)$.

\begin{definition}
Let $\us{n}\in\Fin_\ast$. The stratification $\pi^n_M:M^n_\ast\rightarrow \prt{\us{n}_\ast}$ is defined by sending a point $(x:\us{n}_\ast\rightarrow M_\ast)$ to the unique equivalence relation $R$ providing a factorization
\[
\begin{tikzcd}
\us{n}_\ast\arrow[rr,"x"]\arrow[dr,"f"'] & & M_\ast \\
& \us{n}_\ast/R\arrow[ur,"x'"',hook] &
\end{tikzcd}.
\]
\end{definition}

The proof of \Cref{functoriality_pi} adapts verbatim to prove that the stratification $\pi^n_M:M^n_\ast\rightarrow \prt{\us{n}_\ast}$ is functorial in $\us{n}_\ast$. However, showing that this stratification is conical requires more work than in the setting without boundary.

\begin{lemma}\label{product_of_posets}
Let $R\in \prt{\us{n}_\ast}$. Let $\us{m}_\ast$ be the subset of elements $i$ such that $(i,\ast)\in R$ and $\us{l}_\ast$ the subset of elements $i$ such that $(i,\ast)\not\in R$ or $i=\ast$. Then there is an isomorphism of posets
\[
\prt{\us{n}_\ast}_{R\leq}\cong \prt{\us{m}_\ast}_{R_{|\us{m}_\ast}\leq}\times \prt{\us{l}_\ast}_{R_{|\us{l}_\ast}\leq}.
\] 
\end{lemma}

\begin{proof}
There is a restriction map 
\[
\prt{\us{n}_\ast}_{R\leq}\rightarrow \prt{\us{m}_\ast}_{R_{|\us{m}_\ast}\leq}\times \prt{\us{l}_\ast}_{R_{|\us{l}_\ast}\leq}, (R'\mapsto (R'_{|\us{m}_\ast},R'_{|\us{l}_\ast})).
\]
The converse map 
\[
\prt{\us{m}_\ast}_{R_{|\us{m}_\ast}\leq}\times \prt{\us{l}_\ast}_{R_{|\us{l}_\ast}\leq}\rightarrow \prt{\us{n}_\ast}_{R\leq}
\]
sends a pair $(S,T)$ to the equivalence relation $R':=S\cup T$. More explicitly, $(i,j)\in R'$ if and only if $i,j\in\us{m}_\ast$ and $(i,j)\in S$ or $i,j\in\us{l}_\ast$ and $(i,j)\in T$. A quick check shows that these two maps are inverse of each other and the second one is well defined because whenever $S\subset R_{|\us{m}_\ast},T\subset R_{|\us{l}_\ast}$, we have $S\cup T\subset  R_{|\us{m}_\ast}\cup R_{|\us{l}_\ast}=R$.
\end{proof}

\begin{lemma}\label{product_of_cones}

Let $L_1,\cdots,L_n$ be compact topological spaces. Then there is a homeomorphism $h:C(L_1)\times\cdots\times C(L_n)\cong C(L_1\ast\cdots\ast L_n)$ where $\ast$ denotes the join operation.
\end{lemma}

\begin{proof}
Recall first that the join $L_1\ast\cdots\ast L_n$ consists of all tuples $((x_1,t_1),\cdots,(x_n,t_n))\in (L_1\times I)\times\cdots\times (L_n\times I)$ with the constraint that $\sum_{i=1}^n t_i=1$, quotiented by relations of the form $((x_i,t_i))\sim (x'_i,t_i))$ whenever $t_i=0$. 

Since the topological spaces $L_i$ are compact, one may identify the cones $C(L_i)$ as the quotients 
\[
L_i\times [0,+\infty[/((x,0)\sim \ast_i),
\]
$\ast_i$ being the tip of the cone. Then we get a map 
\[
h:C(L_1)\times\cdots C(L_n)\rightarrow C(L_1\ast\cdots\ast L_n)
\]
defined by sending $(\ast_1,\cdots,\ast_n)$ to $\ast$, the tip of the cone, and sending tuples $((x_i,t_i))$ to 

$\lf(\lf(x_i,\frac{t_i}{\sum_{i=1}^n t_i}\ri),\sum_{i=1}^n t_i\ri)$ when at least one $t_i$ is nonzero. By inspection, this map is well defined and its inverse is given by sending $\ast$ to the tuple $(\ast_1,\cdots,\ast_n)$ and tuples $((x_i,t_i),t)$ to $((x_i,tt_i))$, where $t>0$ is the cone coordinate.

\end{proof}

\begin{proposition}\label{conical_stratification_Mn_pointed}
For any $\us{n}_\ast\in\Fin_\ast$, the stratification $\pi^n_M:M^n_\ast\rightarrow \prt{\us{n}_\ast}$ is conical.
\end{proposition}

\begin{proof}
Fix a point $x=(x_1,\cdots,x_n)\in M^n_\ast$ lying in a stratum $R$. Assume first that $x_i\neq\ast$ for all indices $i$. By \Cref{conical_stratification_Mn}, we get a conical neighbourhood $U\times C(V)\hookrightarrow M^n_\ast$ containing $x$, with $V$ being compact.

Next assume that $x=(\ast,\cdots,\ast)$. Taking a collar $L\times [0,1[$ of the boundary in $M$, we obtain an open neighbourhood of the basepoint $\ast\in M_\ast$ homeomorphic to $C(L)$. Taking products, we obtain an open neighbourhood of $x$ in $M^n$ homeomorphic to $C(L)^n$. This admits a stratification $C(L)^n\rightarrow \prt{\us{n}-\ast}$ as follows: fix the tip of the cone to be its basepoint. Then elements of $C(L)^n$ can be identified with pointed maps $\us{n}_\ast\rightarrow C(L)$ and the $\prt{\us{n}_\ast}$-stratification on $C(L)^n$ is the same as the one from $M^n_\ast$. We then obtain the stratification $C(L^{\ast n})\rightarrow \prt{\us{n}_\ast}$ by precomposing with the homeomorphism $h^{-1}:C(L^{\ast n})\xrightarrow{\cong} C(L)^n$ from \Cref{product_of_cones}. Note that the stratification is invariant in the radial direction: given $((x_i,t_i))\in L^{\ast n}$ and $t,t'>0$, we have $h^{-1}((x_i,t_i),t)=((x_i,tt_i))$ and $h^{-1}((x_i,t_i),t')=((x_i,t't_i))$ which are in the same stratum because all $x_i$'s components are equal and $tt_i=0$ if and only if $t't_i=0$. Moreover, the only point being sent to the minimal element is $x$ itself. This provides the desired conical neighborhood $C(L^{\ast n})\rightarrow \prt{\us{n}_\ast}$ around $x$. 

Finally take $x=(x_1,\cdots,x_n)\in M^n_\ast$ an arbitrary point and $R$ its stratum. Without loss of generality, we may assume that for some index $0\leq m\leq n, x_1=\cdots=x_i=\ast$ and $x_{i+1},\cdots,x_n\neq\ast$. By the two previous cases, we get conical neighbourhoods $C(L^{\ast m})\hookrightarrow \prt{\us{m}_\ast}_{R_{|\us{m}_\ast}\leq}$ around $(x_1,\cdots,x_m)$ in $M^m_\ast$  and $U\times C(V)\hookrightarrow \prt{\us{n-m}_\ast}_{R_{|\us{n-m}_\ast}\leq}$ around $(x_{m+1},\cdots,x_n)$ in $M^{n-m}_\ast$. Taking the product gives a stratified open neighborhood 
\[
U\times C(V)\times C(L^{\ast m})\rightarrow \lf(\prt{\us{m}_\ast}_{R_{|\us{m}_\ast}\leq}\ri)\times \lf(\prt{\us{n-m}_\ast}_{R_{|\us{n-m}_\ast}\leq}\ri)
\]
around $x$ in $M^n_\ast$. Because of \Cref{product_of_posets} and \Cref{product_of_cones}, this amounts to a stratified conical neighborhood $U\times C(V\ast L^{\ast m})\rightarrow \prt{\us{n}_\ast}_{R\leq}$.

\end{proof}

By \Cref{conical_stratification_Mn_pointed} and the functoriality of $\pi^-_M$, one gets a functor 
\[
\text{Enter}^{\prt{\us{-}_\ast}}(M_\ast^-):\Fin_\ast^{op}\rightarrow\Cat^\delta.
\] 

\begin{definition}
We define $(\X_M)^\partial\rightarrow\Fin_\ast$ to be the weighted nerve of the functor \[
\text{Enter}^{\prt{\us{-}_\ast}}(M_\ast^-):\Fin_\ast^{op}\rightarrow\Cat^\delta.
\] 
We define $\con^\partial(M)\subset(\X_M)^\partial$ to be the full subcategory of objects of the form $(\us{n}_\ast,c:\us{n}_\ast\hookrightarrow M_\ast)$.
\end{definition}

With this definition of the configuration category, the construction of the functor $\gamma:\con^<(I)\rightarrow\con^\partial(M)$ is the same as the setting without boundary. Namely, fix a monomorphism of simplicial sets $A\hookrightarrow X$ and a homeomorphism $f:(\norm{A}\hookrightarrow\norm{X})\cong(\partial M\hookrightarrow M)$. This gives a pointed simplicial set $Q:=X/A:\Delta^{op}\rightarrow\Fin_\ast$ and a (noncommutative) diagram

\[
\begin{tikzcd}
\nabla^{op}\arrow[rr,"\text{Enter}^{\pos{\bullet}}{(\Delta^\bullet})"]\arrow[d,"Q"'] & & \Cat \\
\Fin_\ast^{op}\arrow[urr,"\enter{\pos{Q^-}}{\norm{Q^-}}"']
\end{tikzcd}.
\]

Following the same construction as the case without boundary, we get a natural transformation $\enter{\pos{\bullet}}{\Delta^\bullet}\Rightarrow\enter{\pos{Q^-}}{\norm{Q^-}}\circ Q$ and then get the functor $\gamma:\con^<(I)\rightarrow\con^\partial(M)$. The proof of the cofinality is expected to be the same as the one for manifolds without boundary.

The only missing ingredient to recover \Cref{main_simplicial_model} in the boundary setting, is the equivalence with the $\infty$-category of disks embedded in $M$. More formally, following \cite[Section 6.1]{confcat}, let $L$ be a $(d-1)$-dimensional compact manifold without boundary and let $\cMan^\partial$ be the $\infty$-category whose objects are smooth $d$-dimensional manifolds with boundary diffeomorphic to $L$, and for such manifolds $M,N$, the mapping space $\cMan^\partial(M,N)$ is the space of embeddings $M\hookrightarrow N$ fixing the boundary, up to the prescribed diffeomorphisms $\partial M\cong L\cong\partial N$. This admits a full subcategory $\cDisk^\partial\subset\cMan^\partial$ whose objects are manifolds of the form
\[
(L\times [0,1[)\sqcup\bigsqcup\limits_{n}\R^d
\]
for some $n\geq 0$. Then one gets the overcategory $\cDisk^\partial_{/M}$ for an object $M\in\cMan^\partial$. One then would like to establish an equivalence in $(\Cat)_{/\Fin_\ast}$
\[
\con^\partial(M)\simeq\cDisk^{\partial}_{/M}.
\]

Upon proving this equivalence, one would deduce a version of \Cref{main_simplicial_model} for compact manifolds with boundary.

\bibliographystyle{alpha}
\bibliography{bibliography_thesis}

\bigskip
\noindent
\textit{Email address:}
\href{mailto:rodrigue.haya@uclouvain.be}{\texttt{rodrigue.haya@uclouvain.be}}\\
\textsc{Institut de Recherche en Mathématique et Physique,
Université catholique de Louvain, chemin du cyclotron 2 bte L7.01.02, B–1348 Louvain-la-Neuve, Belgium}

\end{document}